\documentclass[12pt]{amsart}
\usepackage{graphicx}
\usepackage{times}
\usepackage{mathptmx} 
\usepackage{amssymb} 
\usepackage{eucal}
\newtheorem{thm}{Theorem}[section]  
\newtheorem{cor}[thm]{Corollary}

\newtheorem{defin}[thm]{Definition} 
\newtheorem{lemma}[thm]{Lemma} 
 
\newtheorem{prop}[thm]{Proposition} 
\newtheorem{conj}{Conjecture} 

\newcommand{\map}{\mbox{$\rightarrow$}}

\newcommand{\eee}{\mbox{$\epsilon$}} 
 
\newcommand{\Ggg}{\mbox{$\Gamma$}}

\newcommand{\bdd}{\mbox{$\partial$}}

\begin{document}  

\title{Generalized Property R and the Schoenflies Conjecture}   

\author{Martin Scharlemann}
\address{\hskip-\parindent
        Martin Scharlemann\\
        Mathematics Department\\
        University of California\\
        Santa Barbara, CA USA}
\email{mgscharl@math.ucsb.edu}

\thanks{Research partially supported by a National Science Foundation grant.  Thanks also to Catalonia's Centre Recerca Matem\`atica for their extraordinary hospitality while this work was being completed} 

\date{\today}

    \begin{abstract}  There is a relation between the generalized Property R Conjecture and the Schoenflies Conjecture that suggests a new line of attack on the latter.  The new approach gives a quick proof of the genus $2$ Schoenflies Conjecture and suffices to prove the genus $3$ case, even in the absence of new progress on the generalized Property R Conjecture.	
\end{abstract}

\maketitle

\section{Introduction and preliminaries}  The Schoenflies Conjecture asks whether every PL (or, equivalently, smooth) $3$-sphere in $S^4$ divides the $4$-sphere into two PL balls.  The appeal of the conjecture is at least $3$-fold:  
\begin{itemize} 

\item The topological version (for locally flat embeddings) is known to be true in every dimension.  Both the PL and the smooth versions (when properly phrased, to avoid problems with exotic structures) are known to be true in every {\em other} dimension.

\item If the Schoenflies Conjecture is false, then there is no hope for a PL prime decomposition theorem for $4$-manifolds, for it would imply that there are $4$-manifolds $X$ and $Y$, not themselves $4$-spheres, so that $X \# Y \cong S^4$.

\item The Schoenflies Conjecture is weaker than the still unsolved $4$-dimensional 
PL Poincar\'e Conjecture, and so might be more accessible.

\end{itemize}

Little explicit progress has been made on the Schoenflies Conjecture for several decades, a time which has nonetheless seen rapid progress in our understanding of both $3$- and $4$-dimensional manifolds.  Here we outline how the Schoenflies Conjecture is connected to another important problem on the border between classical $3$- and $4$-dimensional topology, namely the generalized Property R Conjecture.  We show how how at least some of the last two decades of progress in combinatorial $3$-dimensional topology, particularly sutured manifold theory, can be used to extend the proof of the Schoenflies Conjecture from what are called genus $2$ embeddings of $S^3$ in $S^4$ to genus $3$ embeddings.  In some sense this is a pathetic advance, but it has some philosophical interest:  genus $2$ surfaces have long been known to have special properties (eg the hyperelliptic involution) that are not shared by higher genus surfaces.  That this approach works for genus $3$ suggests that the special properties of genus $2$ surfaces are not needed and so are not a barrier to success for arbitrary, higher genus embeddings.  

We work in the PL category throughout.  All manifolds discussed are orientable.

\section{Generalized Property R}

Recall the famous Property R theorem, proven in a somewhat stronger form by David Gabai \cite{Ga2}:

\begin{thm}[Property R] \label{thm:PropR} If $0$-framed surgery on a knot $K \subset S^3$ yields $S^1 \times S^2$ then $K$ is the unknot. 
\end{thm}

It is well-known (indeed it is perhaps the original motivation for the Property R Conjecture) that Property R has an immediate consequence for the handlebody structure of $4$-manifolds:   

\begin{cor}  \label{cor:PropR} Suppose $U^4$ is a homology $4$-sphere and has a handle structure containing exactly one $2$-handle and no $3$-handles.  Then $U$ is the $4$-sphere.
\end{cor}  

\begin{proof}  Since $U$ has a handle structure with no $3$-handles, dually it has a handle structure with no $1$-handles.  In order for $U$ to be connected, this dual handle structure must then have exactly one $0$-handle, so the original handle structure has a single $4$-handle. 

Let $U_- \subset U$ be the union of all $0$- and $1$-handles of $U$ and $M = \bdd U_-$.  $U_-$ can be thought of as the regular neighborhood of a graph or, collapsing a maximal tree in that graph, as the regular neighborhood of a bouquet of circles.  The $4$-dimensional regular neighborhood of a circle in an orientable $4$-manifold is $S^1 \times D^3$, so $U_-$ is  the boundary
connected sum $\natural_{n} (S^{1} \times D^{3})$, some $n \geq 0$.  (Explicitly, the number of summands $n$ is one more than the difference between the number of $1$-handles and $0$-handles, ie $1 - \chi$, where $\chi$ is the Euler characteristic of the graph.)   It follows that $M = \#_n  (S^{1} \times S^{2})$ and, in particular, $H_2(M) \cong \mathbb{Z}^n$.  Now consider the closed complement $U_+$ of $U_-$ in $U$.  Via the dual handle structure, $U_+$ is obtained by attaching a single $2$-handle to $B^4$, so it deformation retracts to a $2$-sphere and, in particular, $H_2(U_+) \cong \mathbb{Z}$.  
  Since $U$ is a homology $4$-sphere and $H_2(U_-) = 0$, it follows from the Mayer-Vietoris sequence 
$$ H_3 (U) = 0 \map H_2 (M) \map H_2 (U_+) \oplus H_2(U_-) \map H_2(U) = 0$$ 
that $\mathbb{Z} \cong H_2(U_+) \cong \mathbb{Z}^n$, so $n = 1$ and $M = S^1 \times S^2$.

 On the other hand, $U_+$, whose handle structure (dual to that from $U$) consists of a $0$-handle and a $2$-handle,  is visibly the trace of surgery on a knot in $S^3$, namely the attaching map of the $2$-handle.  The framing of the surgery is trivial, since the generator of $H_2(U_+)$ is represented by $* \times S^2 \subset S^1 \times S^2 \cong M$ and this class visibly has trivial self-intersection.  Since the result of $0$-framed surgery on the knot is $M = S^1 \times S^2$, the knot is trivial by Property R (Theorem \ref {thm:PropR}) so $U_+$ is simply $S^2 \times D^2$.  

Hence $U$ is the boundary union $S^1 \times D^3 \cup_{\bdd} S^2 \times D^2$.  Of course the same is true of $S^4$, since the closed complement of a neighborhood of the standard $2$-sphere in $S^4$ is $S^1 \times D^3$.  So we see that $U$ can be obtained from $S^4$ by removing the standard $S^1 \times D^3$ and pasting it back in, perhaps differently.   But it is well-known (and is usefully extended to all 4-dimensional handlebodies by Laudenbach and Poenaru \cite{LP}) that any automorphism of $S^1 \times S^2$ extends to an automorphism of $S^1 \times D^3$, so the gluing homeomorphism extends across $S^1 \times D^3$ to give a homeomorphism of $U$ with $S^4$.
\end{proof}  

The generalized Property R conjecture (cf Kirby Problem 1.82) says this:  

\begin{conj}[Generalized Property R] \label{conj:genR} Suppose $L$ is a framed link of $n \geq 1$ components in $S^3$, and surgery on $L$ via the specified framing yields $\#_{n} (S^{1} \times S^{2})$.  Then there is a sequence of handle slides on $L$ (cf \cite{Ki}) that converts $L$ into a $0$-framed unlink.  
\end{conj}

In the case $n = 1$ no slides are possible, so Conjecture \ref{conj:genR} does indeed directly generalize Theorem \ref{thm:PropR}.  On the other hand, for $n > 1$ it is certainly necessary to include the possibility of handle slides.  For if one starts with the $0$-framed unlink of $n$-components and does a series of possibly complicated handle-slides, the result will be a possibly complicated framed link $L$ of $n$-components.  The result of doing the specified framed surgery on $L$ will necessarily be the same (cf \cite{Ki}) as for the original unlink, namely $\#_n (S^1 \times S^2)$, but $L$ itself is no longer the unlink.  The example $L$ is still consistent with Conjecture \ref{conj:genR} since simply reversing the sequence of handle slides will convert $L$ back to the framed unlink.  So in some sense Conjecture \ref{conj:genR} is the broadest plausible generalization of Theorem \ref{thm:PropR}.  

The generalized Property R Conjecture naturally leads to a generalized Corollary \ref{cor:PropR}:

\begin{prop}  \label{prop:genR} Suppose Conjecture \ref{conj:genR} is true.  Then any homology $4$-sphere $U$ with a handle structure containing no $3$-handles is $S^4$.
\end{prop}

\begin{proof}  Again focus on the $3$-manifold $M$ that separates $U_-$ (the manifold after the $0$ and $1$-handles are attached) from its closed complement $U_+$ in $U$.  The dual handle structure on $U$ shows that $U_+$ is constructed by attaching some $2$-handles to $B^4$.  On the other hand, the original handle structure shows that $U_-$ is the regular neighborhood of a graph, so, as before for some $n$, $U_- \cong  \natural_{n} (S^{1} \times D^{3})$ and $M \cong \#_{n} (S^{1} \times S^{2})$.  In particular $H_2(M) \cong \mathbb{Z}^n$.  Since $U$ is a homology $4$-sphere and $H_2(U_-) = 0$, it follows as before from the Mayer-Vietoris sequence that $H_2(U_+) \cong \mathbb{Z}^n$.  Hence $U_+$ must be obtained from $B^4$ by attaching exactly $n$ $2$-handles.  Then the generalized property R conjecture would imply that $U_+ \cong \natural_{n} (S^{2} \times D^{2})$.  It is shown in \cite{LP} that any automorphism of $\#_{n} (S^{1} \times S^{2}) = \bdd \natural_{n} (S^{1} \times D^{3})$ extends to an automorphism of $\natural_{n} (S^{1} \times D^{3})$.  (This is not quite stated explicitly in \cite{LP} beyond the observation on p. 342, ``mark that no diffeomorphism of $X^p$ was needed here!").  Hence the only manifold that can be obtained by gluing $U_+$ to $U_-$ along $M$ is $S^4$.  \end{proof}

The Proposition suggests this possibly weaker conjecture:

\begin{conj}[Weak generalized Property R conjecture] \label{conj:wgenR}  Suppose attaching $n$ $2$-handles to $B^4$ gives a manifold $W$ whose boundary is $\#_{n} (S^{1} \times S^{2})$.  Then $W \cong \natural_{n} (S^{2} \times D^{2})$. 
\end{conj}

We have then:

\begin{prop}  \label{prop:wgenR}  The weak generalized Property R conjecture (Conjecture \ref{conj:wgenR}) is equivalent to the conjecture that any homology $4$-sphere $U$ with a handle structure containing no $3$-handles is $S^4$.
\end{prop}

\begin{proof}  The proof of Proposition \ref{prop:genR} really required only Conjecture \ref{conj:wgenR}, so only the converse needs to be proved.

Suppose we know that any homology $4$-sphere with a handle structure containing no $3$-handles is $S^4$.  Suppose $W$ is a $4$-manifold constructed by attaching $n$ $2$-handles to $B^4$ and $\bdd W$ is $\#_{n} (S^{1} \times S^{2})$.  Consider the exact sequence of the pair $(W, \bdd W)$: $$ 0 = H_3(W, \bdd W) \map H_2(\bdd W) \map H_2(W) \map H_2(W, \bdd W) \map H_1(\bdd W) \map H_1(W) = 0.$$  Since the last two non-trivial terms are both $\mathbb{Z}^n$, the inclusion induces an isomorphism of the first two non-trivial terms, $H_2(\bdd W) \map H_2(W)$. Attach $ V = \natural_{n} (S^{1} \times D^{3})$ to $W$ by a homeomorphism of their boundaries and call the result $U$.  (There is an obvious homeomorphism of boundaries, and any other one will give the same $4$-manifold, per \cite{LP}).    Then the Mayer-Vietoris sequence for the pair $(W, V)$shows that $U$ is a homology $4$-sphere hence, under our assumption, $U = S^4$.  

$V \subset U$ is just a regular neighborhood of the wedge of $n$ circles $\Ggg$.  Since $U$ is simply connected, $\Ggg$ is homotopic to a standard (ie planar) wedge of circles in $U$ whose complement is 
 $\natural_{n} (S^{2} \times D^{2})$. In dimension $4$, homotopy of 1-complexes implies isotopy (apply general position to the level-preserving map $\Ggg \times I \map U \times I$) so in fact $W \cong \natural_{n} (S^{2} \times D^{2})$ as required. \end{proof}  
 
Setting aside conjecture, here is a concrete extension of Property R:
 
 \begin{prop} \label{prop:Rextend}  Suppose a $2$-handle is attached to a genus $n$ $4$-dimensional handlebody $N = \natural_{n} (S^{1} \times D^{3})$ and the resulting $4$-manifold $N_-$ has boundary  $\#_{n-1} (S^{1} \times S^{2})$.   Then $N_-  \cong \natural_{n-1} (S^{1} \times D^{3})$.
 \end{prop}
 
\begin{proof}  The proof is by induction on $n$; when $n = 1$ this is Property R.  Suppose then that $n > 1$ and let $K \subset \bdd (\natural_{n} (S^{1} \times D^{3})) \cong \#_{n} (S^{1} \times S^{2})$ be the attaching map for the $4$-dimensional $2$-handle.  The hypothesis is then that surgery on $K$ yields $\#_{n-1} (S^{1} \times S^{2})$, a reducible manifold.  But examining the possibilities in \cite{Sch2} we see that this is possible only if $\#_{n} (S^{1} \times S^{2}) - K$ is itself reducible, so in particular one of the non-separating $2$-spheres $\{ * \} \times S^2$ is disjoint from $K$.  Following \cite{LP}, this $2$-sphere bounds a $3$-ball in $N$.  Split $N$ along this $3$-ball, converting $N$ to $\natural_{n-1} (S^{1} \times D^{3})$ and $\bdd N_-$ to $\#_{n-2} (S^{1} \times S^{2})$.  By inductive hypothesis, the split open $N_-$ is $\natural_{n-2} (S^{1} \times D^{3})$ so originally $N_- \cong \natural_{n-1} (S^{1} \times D^{3})$.     
\end{proof}

{\bf Remark:}  Experts will note that, rather than use \cite{Sch2}, one can substitute the somewhat simpler  \cite{Ga1}:  If $n > 1$ then $H_2 (\#_{n} (S^{1} \times S^{2}) - \eta(K)) \neq 0$.  Since both $\infty$- and $0$-framed surgery on $K$ yield reducible (hence non-taut) $3$-manifolds, from \cite{Ga1} it follows that $\#_{n} (S^{1} \times S^{2}) - \eta(K)$ is itself not taut, hence is reducible.
 
 \section{Application:  Heegaard unions}
 
 Let $H^n = \natural_n (S^1 \times D^2)$ denote a $3$-dimensional genus $n$ orientable handlebody and $J^n = \natural_n (S^1 \times D^3)$ denote a $4$-dimensional genus $n$ orientable handlebody. $H^n$ and $J^n$ can also be thought of as regular neighborhoods  in, respectively, $R^3$ and $R^4$ of any graph $\Ggg$ with Euler characteristic $\chi(\Ggg) = 1-n$. 
 
\begin{figure}[tbh]
  \centering
  \includegraphics[width=0.7\textwidth]{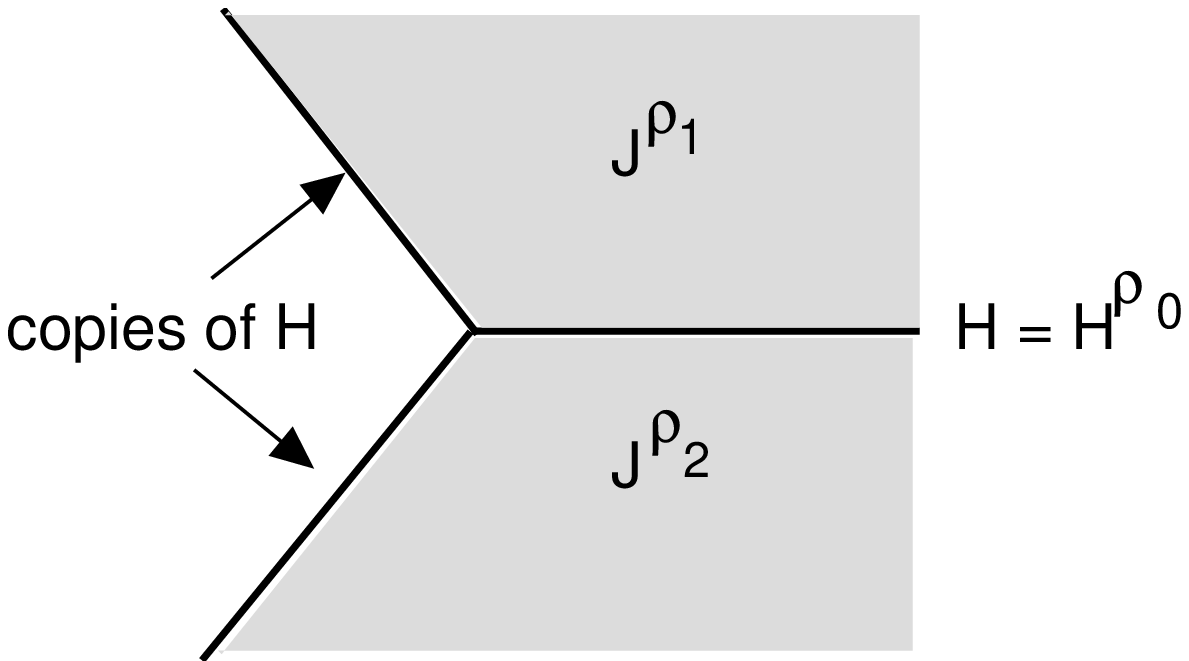}
  \caption{} \label{fig:Heegunion}
  \end{figure}

 \begin{defin}  \label{defin:HeegU}  Suppose, for some $\rho_0, \rho_1, \rho_2 \in \mathbb{N}$,  $H^{\rho_0}$ is embedded into both $\bdd
	J^{\rho_1}$ and $\bdd J^{\rho_2}$ so that its complement in each $\bdd J^{\rho_i}, i = 1, 2$
	is also a handlebody.  Then the
	$4$-manifold $W = J^{\rho_1} \cup_{H^{\rho_0}} J^{\rho_2}$ is called the {\em
	Heegaard union} of the $J^{\rho_i}$ along $H^{\rho_0}$.  See Figure \ref{fig:Heegunion}. 
\end{defin}

  The term Heegaard union comes from the fact that $H^{\rho_0}$ is half of a Heegaard splitting of both $\bdd J^{\rho_1}$ and $\bdd J^{\rho_2}$.
Moreover, if $W$ is such a Heegaard union,  then $(\bdd J^{\rho_1} - H^{\rho_0}) \cup_{\bdd H^{\rho_0}} (\bdd J^{\rho_2} - H^{\rho_0})$ is a Heegaard splitting of $\bdd W$.  The construction here is tangentially related to the construction in \cite[2.4]{BC} of a $4$-dimensional cobordism between three Heegaard-split $3$-manifolds.  Indeed, if two of the three $3$-manifolds in the Birman-Craggs construction are of the form $\#_i (S^1 \times S^2)$ and are then filled in with copies of $\#_i (S^1 \times D^3)$ the result is a Heegaard union.  

\begin{lemma} \label{lemma:rhosum} If a Heegaard union $W = J^{\rho_1} \cup J^{\rho_2}$ is a rational
    homology ball, then $\rho_{0} = \rho_{1} +  \rho_{2}$.  
    \end{lemma}
    
   \begin{proof}  The first and second homology groups (rational
       coefficients) of $W$ are trivial, so the result follows from the Mayer-Vietoris sequence of $W = J^{\rho_1} \cup_{H^{\rho_0}} J^{\rho_2}$:
       $$ H_2 (W) = 0 \map H_1 (H^{\rho_0}) \map H_1 (J^{\rho_1} ) \oplus H_1 (J^{\rho_2}) \map H_1(W) = 0.  $$
\end{proof}
 
 \begin{prop}  \label{prop:HeegU} Suppose a Heegaard union $W = J^{\rho_1} \cup_{H^{\rho_0}} J^{\rho_2}$ is a homology ball and $\bdd W \cong S^3$.  If the weak generalized property R conjecture (Conjecture \ref{conj:wgenR}) is true for $min \{\rho_1, \rho_2 \}$ components, then $W = B^{4}$.
\end{prop}
 
 \begin{proof}  Suppose with no loss of generality that $\rho_1 \leq \rho_2$.  Let $J_i$ denote $J^{\rho_i}, i = 1, 2$ and $H_0$ denote $H^{\rho_0}$.  Consider the genus $\rho_0$ Heegaard splitting of $\bdd J_2$ given by $H_0 \cup_{\bdd H_0} (\bdd J_2 - H_0)$.  According to Waldhausen \cite{Wa} there is only one such Heegaard splitting of $\bdd J_2$ up to homeomorphism, obtained as follows:  Regard $J_2$ as the product of the interval with a genus $\rho _2$ $3$-dimensional handlebody $H$.  Then $H \times \{ 0 \} \subset \bdd (H \times I) = \bdd J_2$ and $\bdd J_2 - (H \times \{ 0 \})$ are both $3$-dimensional handlebodies.  The resulting Heegaard splitting of $\bdd J_2$ is called the {\em product splitting}.  It can be regarded as the natural Heegaard splitting of $\bdd J_2 \cong \#_{\rho_2} (S^{1} \times S^{2})$.  Any other Heegaard splitting (eg the genus $\rho_0$ splitting at hand) is homeomorphic to a stabilization of this standard splitting.  
 
 As proven in \cite{LP} and noted above, any automorphism of $\bdd J_2$ extends over $J_2$ itself, so we may as well assume that the Heegaard splitting $H_0 \cup_{\bdd H_0} (\bdd J_2 - H_0)$ actually is a stabilization of the product splitting.  In particular, and most dramatically, if $\rho_0 = \rho_2$ then no stabilization is required, so $J_2$ is just $H_0 \times I$ and $W \cong J_1$.  Much the same is true if $\rho_0 = \rho_2 + 1$: most of $H_0$ is just $H$, so its attachment to $J_1$ has no effect on the topology of $J_1$.  The single stabilization changes the picture slightly, and is best conveyed by considering what the effect would be of attaching a $4$-ball to $J_1$ {\em not} along one side of the minimal genus splitting of $\bdd B^4$ (ie along $B^3 \subset S^3$), which clearly leaves $J_1$ unchanged, but rather along one side of the once-stabilized splitting of $\bdd B^4$.  That is, $B^4$ is attached to $J_1$ along a solid torus, unknotted in $\bdd B^4$.  But this is exactly a description of attaching a $2$-handle to $J_1$.  So $W$ can be viewed as $J_1$ with a single $2$-handle attached.  In the general situation, in which the product splitting is stabilized $\rho_0 - \rho_2$ times, $W$ is homeomorphic to $J_1$ with $\rho_0 - \rho_2$ $2$-handles attached.  The result now follows from Lemma \ref{lemma:rhosum} and Proposition \ref{prop:wgenR}. \end{proof}
 
 {\bf Remark:}  The link along which the $2$-handles are attached has $\rho_1$ components and, viewed in $S^3$, is part of a genus $\rho_0$ Heegaard splitting.  So its tunnel number can be calculated: $\rho_1 - 1$ tunnels are needed to connect the link into a genus $\rho_1$ handlebody, and another $\rho_0 - \rho_1$ are needed to make it half of a Heegaard splitting.  Hence the tunnel number is $\rho_0 - 1$.  This fact may be useful, but anyway explains why \cite{Sch1} could be done just knowing Property R for tunnel number one knots.
 
  \begin{cor}  \label{cor:HeegU} Suppose a Heegaard union $W = J^{\rho_1} \cup_{H^{\rho_0}} J^{\rho_2}$ is a homology ball and $\bdd W \cong S^3$.  If $\rho_0 \leq 3$ then $W = B^{4}$.
\end{cor}

\begin{proof}  By Lemma \ref{lemma:rhosum}, $\rho_1 + \rho_2 \leq 3$, hence $min \{\rho_1, \rho_2 \} \leq 1$.  The result then follows from Proposition \ref{prop:HeegU} and Theorem \ref{thm:PropR}. \end{proof}

 \section{Handlebody structure on 3-manifold complements}
 
Suppose $M \subset S^4$ is a connected closed PL or smooth $3$-submanifold.  In this section we discuss the handlebody structure of each complementary component of $M$.  
 
 It is a classical result (cf \cite{KL}) that $M$ can be isotoped so that it is in the form of a {\em rectified critical level embedding}.  We briefly review what that means.  
 
 Informally, the embedding $M \subset S^4$ is in the form of a critical level embedding if it has a handlebody structure in which each handle is horizontal with respect to the natural height function on $S^4$, and $M$ intersects each region of $S^4$ between handle levels in a vertical collar of the boundary of the part of $M$ that lies below (or, symmetrically, above).  More formally, regard $S^4$ as the boundary of $D^4 \times [-1, 1]$, so $S^4$ consists of two $4$-balls $D^4 \times \pm 1$ (called {the poles}) added to the ends of  $S^3 \times [-1, 1]$.  Let $p: S^3 \times [-1, 1] \map [-1, 1]$ be the natural projection.  For $-1 < t < 1$ denote $p^{-1}(t)$ by $S^3_t$.  Then $M \subset S^3 \times [-1, 1] \subset S^4$ is a critical level embedding if there are a collection $t_1 < t_2 < \ldots < t_n$ of values in $(-1, 1)$ and a collection of closed surfaces $F_1, \ldots F_n \subset S^3$ so that 
 
 \begin{enumerate}
 
 \item $p(M) = [t_1, t_n] \subset (-1, 1)$
 
 \item for each $1 \leq i \leq n-1$, $M \cap (S^3 \times (t_i, t_{i+1})) = F_i \times (t_i, t_{i+1})$
 
 \item $M \cap S^3_{t_1} = B^3$ with boundary $F_1$
 
 \item For each $2 \leq i \leq n$, $F_{i}$ is obtained from $F_{i - 1}$ by a $j$-surgery, some $0 \leq j \leq 3$.  That is, there is a $3$-ball $D^j \times D^{3 - j} \subset S^3$ incident to $F_{i-1}$ in $\bdd D^j \times D^{3 - j}$ and  $F_i$ is obtained from $ F_{i - 1}$ by replacing $\bdd D^j \times D^{3 - j}$ with $D^j \times \bdd D^{3 - j}$.

  \item for each $2 \leq i \leq n$, $M \cap S^3_{t_i}$ is the trace of the surgery above.  That is, it is the union of $F_{i - 1}$, $F_i$ and $D^j \times D^{3 - j}$.   
 \end{enumerate}
 
\begin{figure}[tbh]
  \centering
  \includegraphics[width=0.9\textwidth]{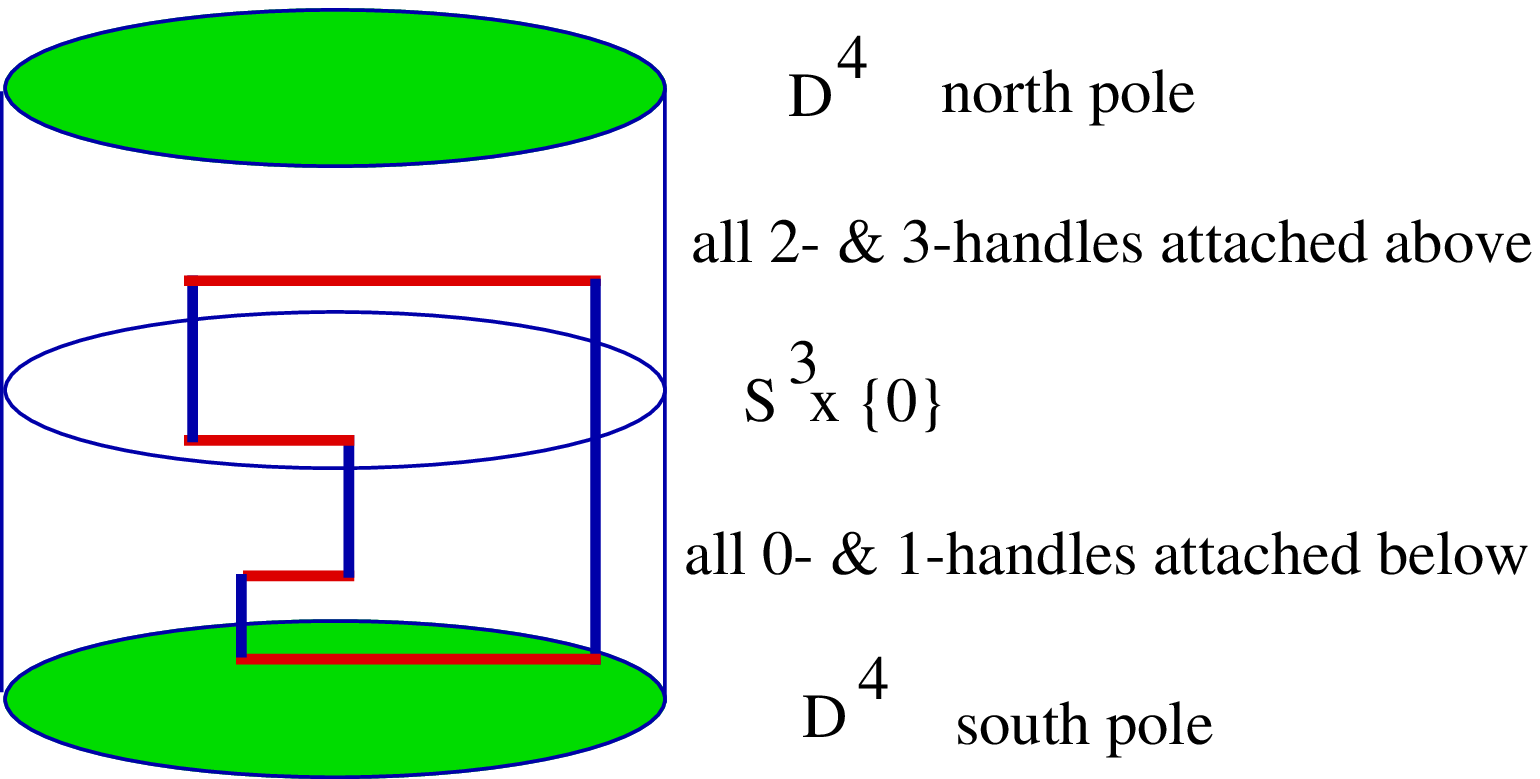}
  \caption{} \label{fig:rectified}
  \end{figure}

Such an embedding gives rise to a handle structure on $M$ with $n$ handles added successively at levels $t_1, \ldots , t_n$.  $j$ is the index of the handle $D^j \times D^{3 - j}$.  A critical level embedding is called {\em rectified} if, for  $0 \leq j \leq 2$, each handle of index $j$ occurs at a lower level than each handle of index $j+1$.  Furthermore, all $0$- and $1$-handles lie below $S^3_0$ and all $2$- and $3$-handles lie above $S^3_0$.  See Figure \ref{fig:rectified}.

Note that the surface $M \cap S^3_0$ is a Heegaard surface for $M$, since all $0$- and $1$-handles lie on one side (namely in $S^3 \times [-1, 0]$) and all $2$- and $3$-handles lie on the other ($S^3 \times [0, 1]$). In particular,  $M \cap S^3_0$ is connected.  It is easy to see, \cite[Lemma 1.4]{Sch1}, though not completely obvious, that if the first $1$-handle attached to the boundary of a $0$-handle is incident to the $0$-handle at only one end, then the handles cancel and there is a rectified embedding of $M$ in which neither handle appears.  So, minimizing the number of handles, we will henceforth assume that the first $1$-handle incident to each $0$-handle is incident to it in both ends.  Equally important is the dual to this remark:  the boundary of the core of any $2$-handle is essential in the surface to which the $2$-handle is attached. To summarize:

\begin{lemma} \label{lemma:convention}  Any rectified critical level embedding of $M$ may be isotoped rel $M \cap S^3_0$ to a rectified critical level embedding with no more (but perhaps fewer) handles of any index, such that 
\begin{itemize}
\item the first $1$-handle incident to each $0$-handle is incident to it in both ends and
\item the core of any $2$-handle attached in $S^3_t$ is a compressing disk for $M \cap S^3_{t - \eee}$.
\end{itemize}
\end{lemma}

We will henceforth consider only rectified critical level embeddings with these two properties.  

\begin{defin}  The {\em genus} of the embedding of $M$ in $S^4$ is the genus of the Heegaard surface $M \cap S^3_0$.
\end{defin}

\bigskip

It will be important to understand how a rectified critical level embedding induces a handlebody structure on each of its closed complementary components $X $ and $Y$.  Let $X$ denote the component of $S^4 - M$ that contains the south pole $D^4 \times \{ -1 \}$.  For each generic $t \in (-1, 1)$ let $Y^-_t$, (resp $X^-_t $, $M^-_t$) be the part of $Y$ (resp $X$, $M$) lying below level $t$ or, more formally, the $4$-manifold with boundary $Y \cap (S^3 \times [-1, t])$ (resp $X \cap (S^3 \times [-1, t])$, $3$-manifold with boundary $M \cap (S^3 \times [-1, t]$). Symmetrically, let $Y^+_t$, (resp $X^+_t $, $M^+_t$) be the part of $Y$ (resp $X$, $M$) lying above level $t$, that is, the $4$-manifold with boundary $Y \cap (S^3 \times [ t, 1])$ (resp $X \cap (S^3 \times [ t, 1])$, $3$-manifold with boundary $M \cap (S^3 \times [ t, 1])$).  Finally, let $Y^*_t$, (resp $X^*_t $, $M^*_t$) be the part of $Y$ (resp $X$, $M$) lying at level $t$, that is the $3$-manifold with boundary $Y \cap S^3_t$ (resp $X \cap S^3_t$, closed surface $M \cap S^3_t$). Thus $\bdd Y^-_t$ is the union of $M^-_t$ and $Y^*_t$.  If $t_i < t < t_{i+1}$ then $\bdd M^-_t =  \bdd M^+_t  = M^*_t  = F_i \subset S^3_t$ and $Y^*_t$ consists of a collection of closed complementary components of $F_i$ in $S^3_t$.  Each component of $F_i$ in $S^3_t$ is incident to $Y^*_t$ on exactly one side and to $X^*_t$ on the other.

Clearly as long as no $t_i$ lies between the values $t < t'$, then $Y^{\pm}_t \cong Y^{\pm}_{t'}$, since the region between them is just a collar on part of the boundary.  On the other hand, for each $t_i$, consider the relation between  $Y^{-}_{t_i - \eee}$ and $Y^{-}_{t_i + \eee}$.  We know that $F_i$ is obtained from $F_{i-1}$ by doing $j$-surgery along a $j$-disk in $S^3 - F_{i-1}$.  If that $j$-handle lies on the $Y$ side of $F_{i-1}$ (in the sphere $S^3_{t_i - \eee}$) then $Y^*_{t_i + \eee}$ is homeomorphic to just $Y^*_{t_i - \eee}$ with that $j$-handle removed.  So $Y^-_{t_i + \eee}$ is still just $Y^-_{t_i - \eee}$ with a collar added to part of its boundary, but only to the complement of the $j$-handle in $Y^*_{t_i - \eee}$.  Hence it is still true that $Y^-_{t_i + \eee} \cong Y^-_{t_i - \eee}$.  On the other hand, if the $j$-handle lies on the $X$ side of $F_{i-1}$, then $Y^-_{t_i + \eee}$ is homeomorphic to $Y^-_{t_i - \eee}$ but with a (4-dimensional) $j$-handle added, namely the product of the interval $[t_i, t_i + \eee]$ with the $3$-dimensional $j$-handle added to $M^-_{t_i - \eee}$ in $S^3_{t_i}$.  

\begin{figure}[tbh]
  \centering
  \includegraphics[width=0.6\textwidth]{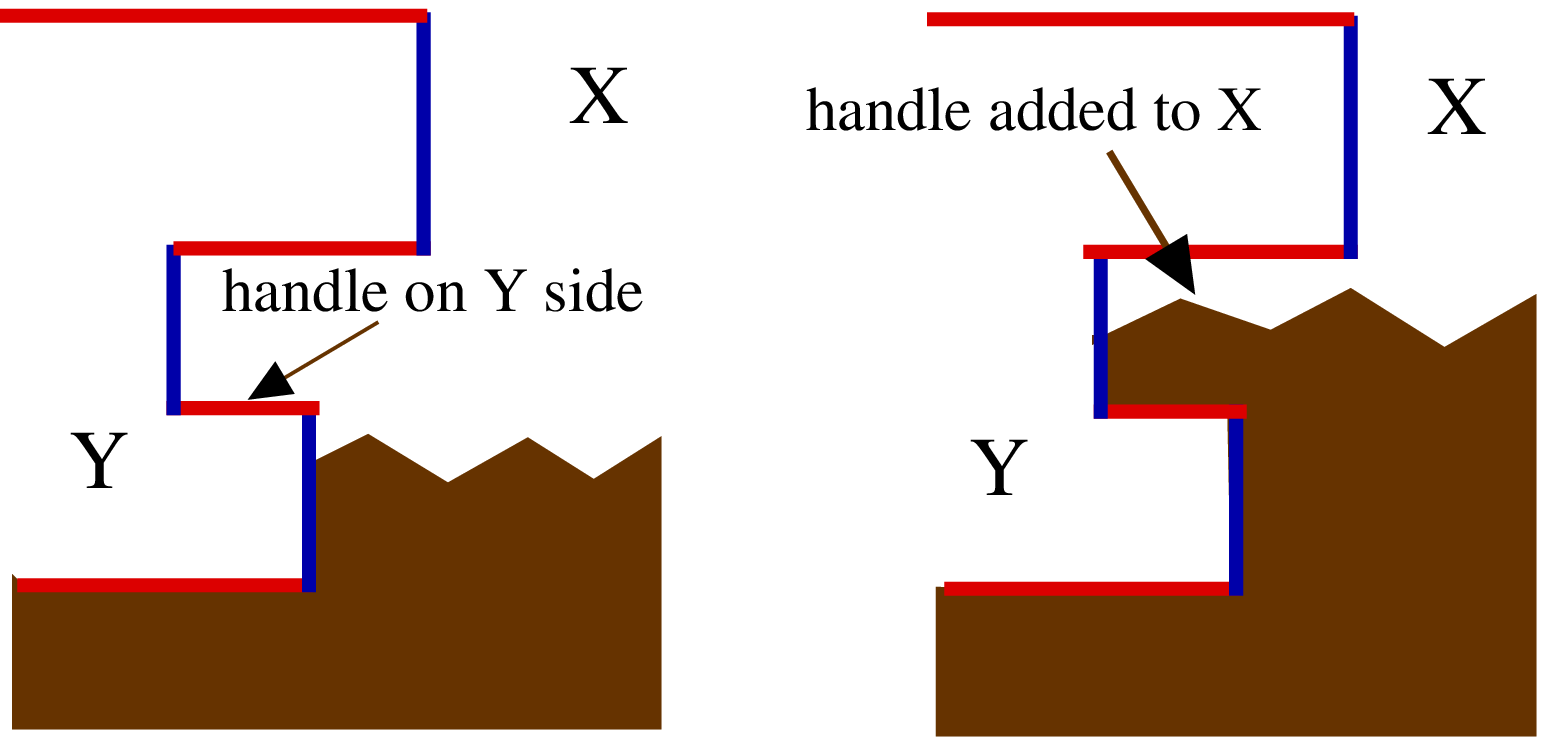}
  \caption{} \label{fig:risingwater}
  \end{figure}

We have then the general rule, sometimes called the rising water rule (cf Figure \ref{fig:risingwater}:

\begin{lemma} \label{lemma:handlerule} 
\begin{enumerate}

\item If the $j$-surgery at level $t_i$ has its core in $Y$, then $Y^-_{t_i + \eee} \cong Y^-_{t_i - \eee}$.

\item If the $j$-surgery at level $t_i$ has its core in $X$, then $Y^-_{t_i + \eee} \cong Y^-_{t_i - \eee}$ with a $j$-handle attached.

\end{enumerate}
\end{lemma}

Of course the symmetric statements hold for $X$.  Note that since $X$ contains the south pole, $X_{t_0 - \eee} \cong B^4$ whereas $Y_{t_0 - \eee} = \emptyset$.  Just as $M^*_0$ is a Heegaard surface for $M$, $X^-_0$ and $Y^-_0$ are connected $4$-manifolds, constructed from just $0$- and $1$-handles.  In other words, there are integers $n_x, n_y \geq 0$ so that  $X_0 = \natural_{n_x} (S^{1} \times D^{3})$ and $Y_0 = \natural_{n_y} (S^{1} \times D^{3})$.  

Each handle in $M^-_0$ corresponds to a handle of the same index in exactly one of $X^-_0$ or $Y^-_0$, so there is a connection between $n_x, n_y$ and the genus $g$ of $M^*_0$:  The critical level embedding defines a handlebody structure on $M^-_0$ with $a$ $0$-handles and $b$ $1$-handles, where $$b - a + 1 = g.$$ If $a > g$ then there would be at least one $0$-handle in the critical level embedding which is first incident to a $1$-handle on a single one of its ends, violating the Handle Cancellation Lemma \ref{lemma:convention}.   So $$ a \leq g.$$  

Let $a_x, a_y$ (resp $b_x, b_y$) denote the number of $0$- (resp $1$-) handles in the critical level embedding whose cores lie in $X$ and $Y$.  We have from above that $a_x + a_y = a, b_x + b_y = b, n_x = b_y - a_y$ and $n_y = b_x - a_x + 1$. (The asymmetry is explained by noting that the south pole is a $0$-handle for $X$.)  It follows that $$n_x + n_y = g.$$  

Another way of counting $n_x$ and $n_y$ is this:  Suppose a $1$-handle at critical level $t_i$ has its core lying in $X$, say.  If the ends of the $1$-handle lie in distinct components of $F_{i-1}$ then the $1$-handle adds a $1$-handle to $Y$ but nothing to its genus.  In contrast, if the ends of the $1$-handle lie on the same component of $F_{i-1}$ then it adds $1$ to the genus of $Y$.  A count of the total number of the latter sort of $1$-handles lying in $X$ (resp $Y$) gives $n_y$ (resp $n_X$).  

\begin{figure}[tbh]
  \centering
  \includegraphics[width=0.6\textwidth]{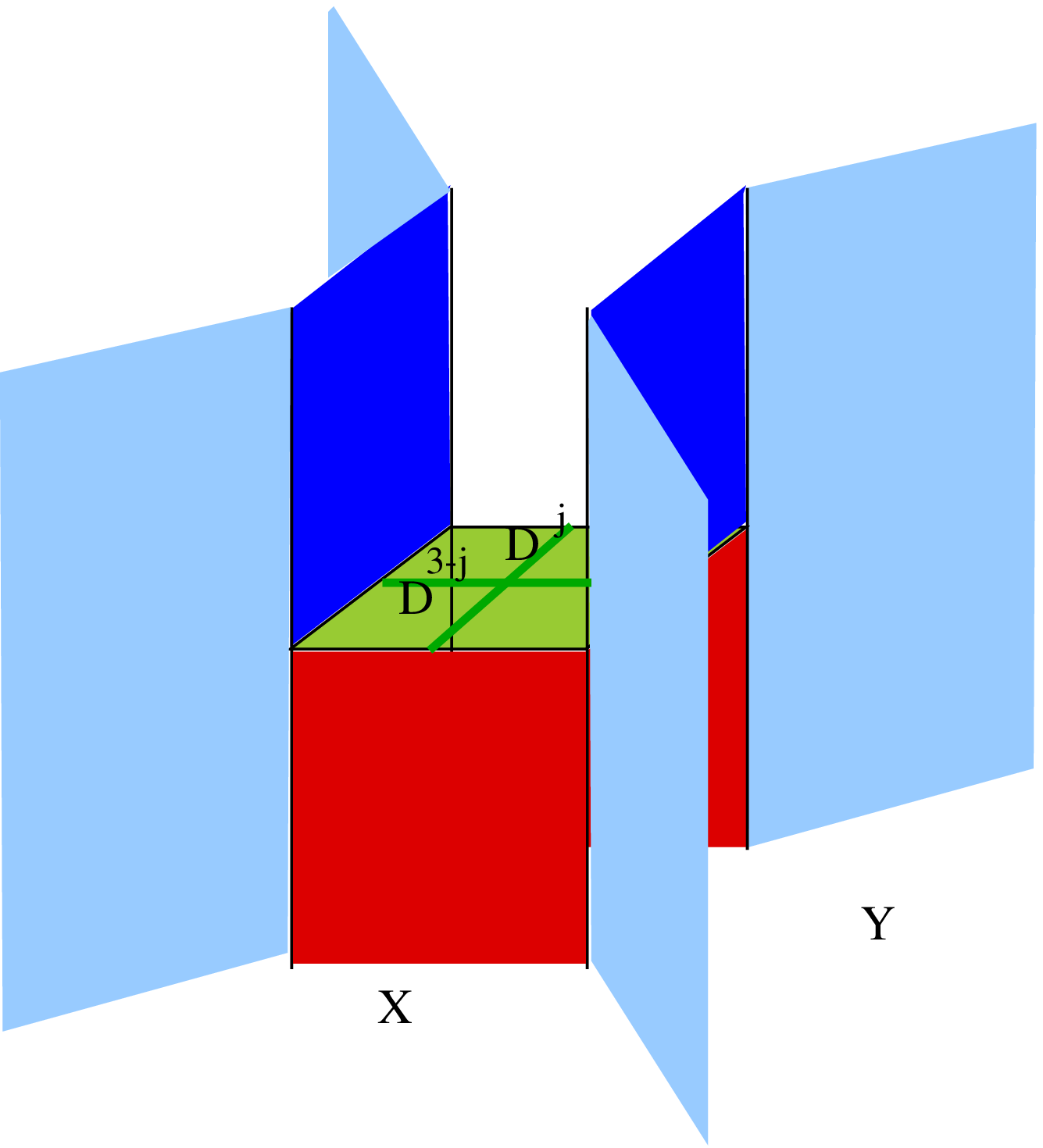}
  \caption{} \label{fig:dualityin3}
  \end{figure}

For everything that has been said about $X^-$ and $Y^-$ there is a dual statement for $X^+$ and $Y^+$, easily obtained by just inverting the height function.  The result is that, beyond the standard $4$-dimensional duality of handle structures on $X$ and on $Y$, there is a kind of $3$-dimensional duality between handles in $X$ and handles in $Y$, induced by the $3$-dimensional duality of handles in $M$.  See Figure \ref{fig:dualityin3}.

To be concrete: $X^+_0$ is also a solid $4$-dimensional handlebody.  To determine its genus, consider the core of each $2$-handle, say at critical height $t_i > 0$.  If the core of the $2$-handle lies on the $X$ side of $F_{i-1}$ then the cocore lies on the $Y$ side of $F_i$ so it corresponds to a $1$-handle in $X^+_0$.  This $1$-handle adds genus to $X^+_0$ if and only if the boundary of the $2$-handle is non-separating in $F_{i-1}$.

\begin{figure}[tbh]
  \centering
  \includegraphics[width=0.6\textwidth]{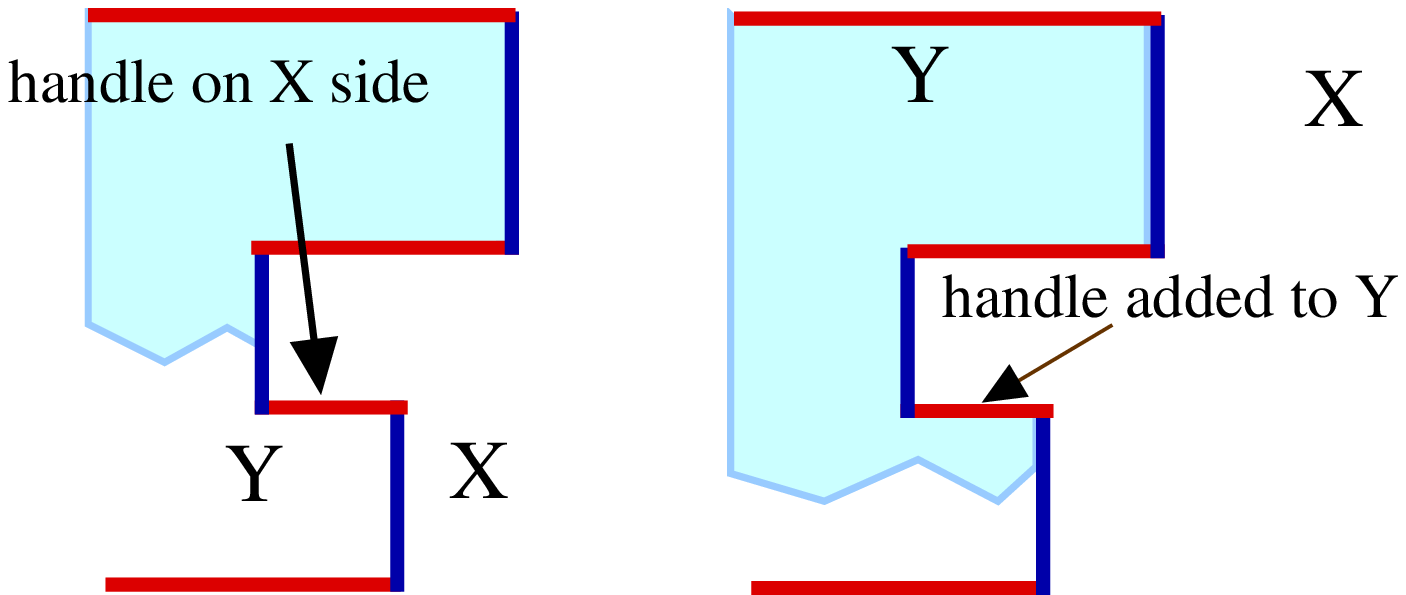}
  \caption{} \label{fig:hydrogen}
  \end{figure}

To see how this occurs, consider the ``dual rule'' to Lemma \ref{lemma:handlerule}.  That is, suppose again that $F_i$ is obtained from $F_{i-1}$ by doing $j$-surgery along a $j$-disk in $S^3 - F_{i-1}$ and ask how $Y^{+}_{t_i - \eee}$ and $Y^{+}_{t_i + \eee}$ differ.  If the $j$-surgery at level $t_i$ has its core in $Y$, then, viewed from above instead of below, there is a corresponding $3-j$ surgery with its core in $X$.  So, following the argument of Lemma \ref{lemma:handlerule}, $Y^+_{t_i - \eee} \cong Y^+_{t_i + \eee}$ with a $(3-j)$-handle attached.  On the other hand, if the core of the $j$-handle lies in $X$, $Y^+$ is unchanged.  This might be called the descending hydrogen rule (cf Figure \ref{fig:hydrogen}). 

To summarize all possibilities:

\begin{lemma} \label{lemma:dualrule} Suppose $F_i$ is obtained from $F_{i-1}$ by doing $j$-surgery along a $j$-disk in $S^3 - F_{i-1}$

\begin{enumerate}
\item If the $j$-surgery at level $t_i$ has its core in $Y$, then 
\begin{itemize}  
\item $Y^-_{t_i + \eee} \cong Y^-_{t_i - \eee}$
\item $X^-_{t_i + \eee} \cong X^-_{t_i - \eee}$ with a $j$-handle attached
\item $Y^+_{t_i - \eee} \cong Y^+_{t_i + \eee}$ with a $(3-j)$-handle attached
\item $X^+_{t_i - \eee} \cong X^+_{t_i + \eee}$.
\end{itemize}

\item If the $j$-surgery at level $t_i$ has its core in $X$, then 
\begin{itemize}
\item $Y^-_{t_i + \eee} \cong Y^-_{t_i - \eee}$ with a $j$-handle attached
\item $X^-_{t_i + \eee} \cong X^-_{t_i - \eee}$
\item $Y^+_{t_i - \eee} \cong Y^+_{t_i + \eee}$
\item $X^+_{t_i - \eee} \cong X^+_{t_i + \eee}$ with a $(3-j)$-handle attached.
\end{itemize}

\end{enumerate}
\end{lemma}

Here is a simple example of how this $3$-dimensional duality can be useful:

\begin{prop} \label{prop:oneside}  Suppose there is a rectified critical level embedding of $M = S^3$ in $S^4$ so that the $0$- and $1$-handles, as they are successively attached, all lie on the $X$-side.  Then $X \cong B^4$.
\end{prop}

\begin{proof} Following Lemma \ref{lemma:dualrule}, $X$ has no $0$ or $1$-handles, so it only has $2$- and $3$-handles.  Dually (in the standard $4$-dimensional handle duality of $X$), $X$ can be constructed with only $1$ and $2$-handles.  Neither of these statements, in itself, is enough to show that $X$ is a $4$-ball.

Consider, however, what the given information tells us about $Y$, following Lemma \ref{lemma:dualrule} applied to the construction of $Y$ from above:  The possible $2$- and $3$-handles in the construction of $X$ from below correspond respectively to $1$- and $0$- handles in the construction of $Y$ from above.  Similarly, the lack of $0$- and $1$-handles (beyond the south pole) for $X$ constructed from below corresponds to a lack of $3$- and $2$-handles for $Y$ when constructed from above.  Hence $Y$ has only $0$- and $1$-handles, ie it is a $4$-dimensional handlebody.  On the other hand, because it is the complement of $S^3$ in $S^4$ it is a homotopy $4$-ball, so the handlebody must be of genus $0$, ie $Y$ is a $4$-ball.  Then its complement $X$ is also a $4$-ball.  \end{proof}

\section{Two proofs of the genus two Schoenflies Conjecture}

Informed by the ideas above, we present two proofs of the genus $2$ Schoenflies Conjecture.  The first is similar in spirit (though different in detail) to the original proof of \cite{Sch1}.  The second uses a  different approach, one that aims to simplify the picture by reimbedding $X$ or $Y$.  

Here is a more general statement, relevant to the classical approach:

\begin{prop} \label{prop:0handle} Suppose a $3$-sphere $M$ has a genus $g$ rectified critical level embedding in $S^4$ with at most two $0$-handles or at most two $3$-handles.   If the generalized property R conjecture is true for links of $g - 1$ components then $M$ divides $S^4$ into two $PL$ $4$-balls.
\end{prop}

\begin{proof} Perhaps inverting the height function, assume without loss of generality that $M$ has at most two $3$-handles.  The roles of $X$ and $Y$ can be interchanged by passing the lowest $0$-handle over the south pole, so we can also assume without loss of generality that the first (that is, the lowest) $3$-handle for $M$ lies in $Y$ and so represents the addition of a $3$-handle to $X$.  The second $3$-handle (and so the last handle) of $M$ either lies in $X$ or in $Y$, but these options are isotopic by passing the handle over the north pole.  So, via an isotopy of this handle, we can choose whether both $3$-handles of $M$ lie in $Y$ (and so represent attaching of $3$-handles to $X$ and not $Y$) or one each lies in $X$ and $Y$.  See Figure \ref{fig:low3handle}.  

\begin{figure}[tbh]
  \centering
  \includegraphics[width=0.8\textwidth]{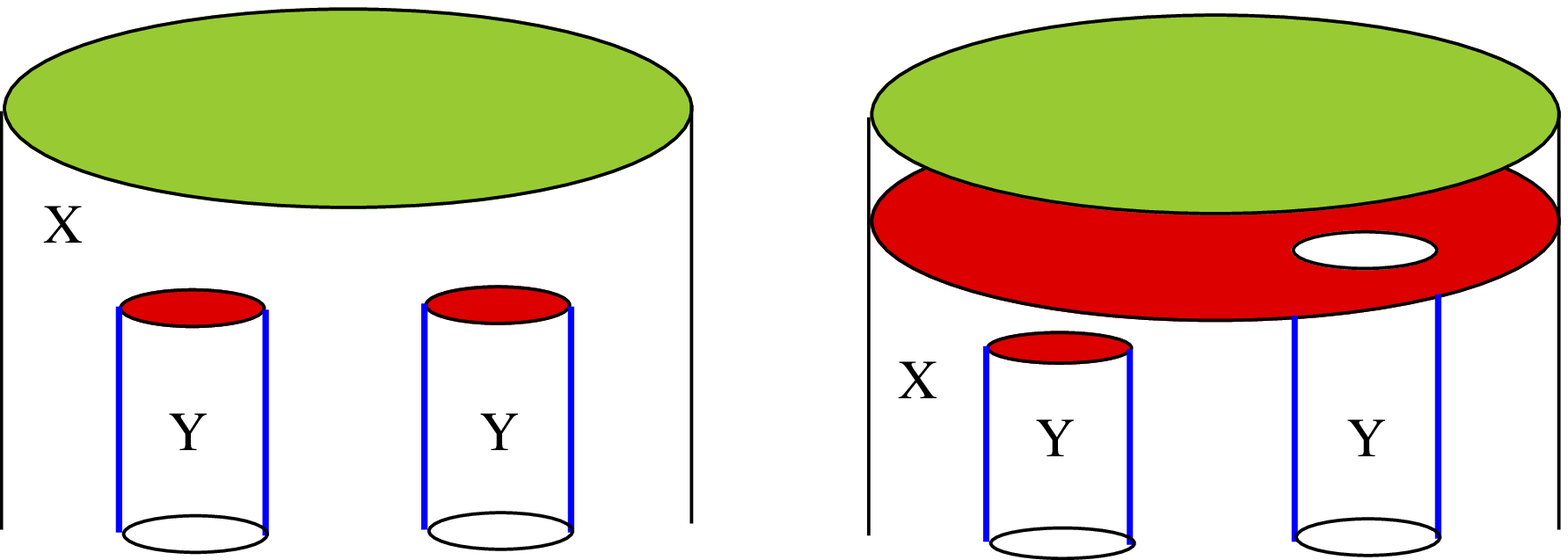}
  \caption{} \label{fig:low3handle}
  \end{figure}  

Now consider the genera $n_x$ and $n_y$ of the $4$-dimensional handlebodies $X^-_0, Y^-_0$, with $g = n_x + n_y$.  If $n_x = 0$ then $X^-_0$ is a $4$-ball.  $X$ is obtained from this $4$-ball by attaching some number of $2$ and $3$-handles, and also a $4$-handle if the north pole of $S^4$ lies in $X$.  There are as many, total, of $2$- and $4$- handles as there are $3$-handles (since $X$ is a homotopy $4$-ball) and the argument of the previous paragraph ensures that we can arrange it so that $X$ contains at most one $3$-handle.  Viewed dually, this means that $X$ can be constructed from $\bdd X = S^3$ with no $3$-handles, and at most one each of $1$- and $2$-handles.  The result then follows from Corollary \ref{cor:PropR}.   

If $n_x \geq 1$ then $n_y \leq g-1$ and, first arranging as above so that $Y$ has no $3$-handles, the result again follows from the proof of Proposition \ref{prop:genR}. \end{proof}

\begin{cor}[Sch]  \label{cor:genus2} Each complementary component of a genus $2$ embedding of $M = S^3$ in $S^4$ is a $4$-ball.
\end{cor}

\begin{proof} As noted above, we can assume that the number $a$ of $0$-handles in the rectified embedding of $M$ is no larger than $g = 2$.  Proposition \ref{prop:0handle} then shows the result follows from Property R, via Corollary  \ref {cor:PropR} . \end{proof}

The reimbedding proof of the genus $2$ Schoenflies Conjecture begins with a more general claim that follows from our results above for Heegaard unions:

\begin{prop} \label{prop:midHeeg} Suppose $M \cong S^3$ has a rectified critical level embedding in $S^4$ so that $Y^*_0$ (resp $X^*_0$) is a \underline{handlebody} of genus $\rho_0$.  If the generalized property R conjecture is true for $[\rho_0/2]$ components then $Y \cong B^4$ (resp $X \cong B^4$).
\end{prop}

\begin{proof}  It was noted above that $Y^-_0$ is a $4$-dimensional handlebody and $M^-_0$ is a $3$-dimensional handlebody.  The latter fact, and the hypothesis, imply that $M^-_0 \cup Y^*_0$ is a Heegaard splitting of $\bdd Y^-_0$.  Viewing the critical level embedding from the top down we symmetrically see that $Y^+_0$ is a $4$-dimensional handlebody and $Y$ is a Heegaard union of $Y^-_0$ and $Y^+_0$ along $Y^*_0$.  

Let $\rho_1, \rho_2$ denote the genera of $Y^-_0$ and $Y^+_0$ respectively.  Since $M$ is a $3$-sphere, each complementary component of $M$ is a homotopy $4$-ball.  In particular, following lemma \ref{lemma:rhosum}, $\rho_1 + \rho_2 = \rho_0$.  The result now follows from Proposition \ref{prop:HeegU}. \end{proof}

Proposition \ref{prop:midHeeg} suggests a clear strategy for a proof of the general Schoenflies Conjecture, assuming the generalized Property R Conjecture:  Given a rectified critical level embedding of $M = S^3$ in $S^4$, try to reimbed $X$ (or $Y$), still a rectified critcal level embedding, so that afterwards, either the $3$-manifold $X^*_0$ or its complement $Y^*_0$ is a handlebody.  Or at least  more closely resembles a handlebody.  For even if a series of reimbeddings, first of $X$, then of its new complement $Y'$, then of the new complement of $Y'$, etc, eventually leads to a handlebody cross-section at height $0$, we are finished.  For once one of the complementary components of the multiply reimbedded $M$ is a $4$-ball, we have that both are, hence the previous complementary components, in succession, leading back to the original $X$ and $Y$ are all $4$-balls.  (This is more formally explained in the proof of Corollary \ref{cor:final}.) What follows is a proof of the genus $2$ Schoenflies Conjecture built on this strategy.

In order to be as flexible as possible in reimbedding $X$ or $Y$ we first prove a technical lemma which roughly shows that, at the expense of some vertical rearrangement of the $3$-handles (or, dually, the $0$-handles), the core of a $2$-handle (resp, the cocore of a $1$-handle) can be moved rel its boundary to another position without affecting the isotopy class of $M$ or even the embedding of $M$ below the specified $2$-handle (resp above the specified $1$-handle).

Suppose, as above, $M$ is a rectified critical level embedded $S^3$ in $S^4$.  

\begin{lemma}[Prairie-Dog Lemma] \label{lemma:prairie}
Let $E \subset S^3_{t_i} - F_{i-1}$ be the core of the $2$-handle added to $F_{i-1}$ at critical level $t_i > 0$ and let $t$ be a generic height such that $t_{i-1} < t < t_i$.  Let $E' \subset S^3_{t_i} - F_{i-1}$ be another disk, with $\bdd E'$ isotopic to $\bdd E$ in $F_{i-1}$.  Then there is a proper isotopy of $M^+_t$ in $S^3 \times [t, 1]$ so that afterwards 
\begin{itemize}
\item the new embedding $M'$ of $S^3$ is still a rectified critical level embedding
\item the critical levels and their indices are the same for $M$ and $M'$
\item the core of the $2$-handle at critical level $t_i$ is $E'$ and
\item for any generic level $t$ below the level of the first $3$-handle, $M^-_t \cong M'^-_t$.
\end{itemize}

\end{lemma}

\begin{proof}  With no loss we take $\bdd E'$ parallel (hence disjoint) from $\bdd E$.  Let  $k$ be the number of $2$-handles above level $t_i$ and $n = |E \cap E'|$.  The proof is by induction on the pair $(k, n)$, lexicographically ordered.

\bigskip

{\bf Case 1:} $k = n = 0$

\bigskip

In this case,  $2$-handle attached at level $t_i$ is the last $2$-handle attached and $E'$ is disjoint from $E$.  Then the union of $E$ and $E'$ (and a collar between their boundaries) is a $2$-sphere in $S^3 - F_{i-1}$; if it bounds a $3$-ball in $S^3 - F_{i-1}$ then the disks are isotopic and there is nothing to prove.  If it does not bound a $3$-ball, let $S$ be the parallel reducing sphere for $S^3 - F_{i}$ and $B \subset S^3$ be the ball it bounds on the side that does not contain the component to which $E$ and $E'$ are attached.  Since $t_i$ is the highest $2$-handle, each component of $F_{i} \cap B$ is a sphere and each is eventually capped off above $t_i$ by a $3$-handle.  

If all are capped off by 3-handles that lie within $B$, push all of $M \cap (B \times [t, 1])$ vertically down to a height just above $t$ so that afterwards, $E$ is isotopic to $E'$ in $S^3_{t_i - \eee} - M$.  Perform the isotopy, then push $M \cap (B \times [t, 1])$ back up, so that the $3$-handles are attached above height $t_i$.  The number of $3$-handles attached (namely, the number of components of $F_i$) is the same, so, although perhaps rearranged in order, the critical heights at which the $3$-handles are attached can be restored to the original set of critical heights.  See Figure \ref{fig:prairiedog}.

 \begin{figure}[tbh]
  \centering
  \includegraphics[width=1.0\textwidth]{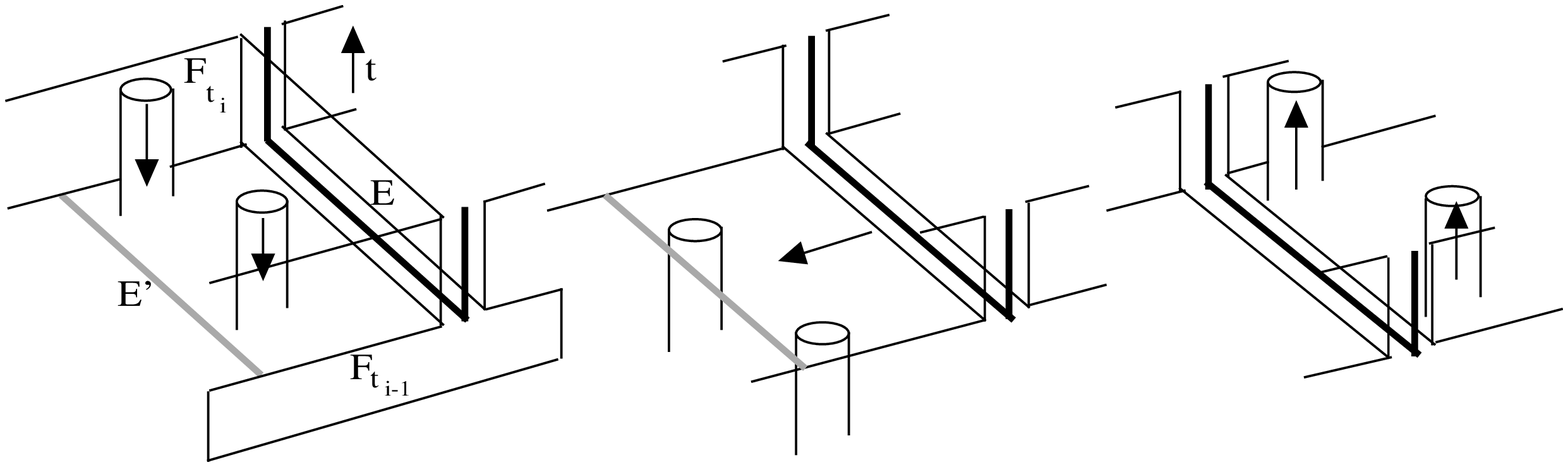}
  \caption{} \label{fig:prairiedog}
  \end{figure}

The picture is only a little changed if one of the components of $F_i$ lying in $B$ is eventually capped off by a $3$-handle not in $B$.  The proof is by induction on the number of such handles.  Consider the highest such handle, say at level $t_j$, capping off a sphere component $S$ of $F_j \cap B$ in $S^3_{t_j}$.  Let $B'$ be the $3$-ball component of $S^3_{t_j} - S$ that does lie in $B$, ie the complement of the $3$-handle in $S^3_{t_j}$.  If there are no components of $F_j$ in the interior of $ B'$ then the $3$-handle is isotopic to $B'$ via passing over the north pole.  This isotopy decreases by one the number of $3$-handles not lying in $B$, completing the inductive step.  If some components of $F_j$ do lie in $B'$ note that eventually they are capped off by $3$-handles lying in $B'$ (by choice of $t_j$).  Simply push $M \cap (B' \times [t_j, 1])$ down below level $t_j$ and do the pole pass described above.  

\bigskip

{\bf Case 2:} $k = 0, n > 0$

\bigskip

Consider an innermost disk $E'_0 \subset E'$ cut off by $E$ in $E'$.  Then the union of $E'_0$ and the subdisk $E_0$ of $E$ bounded by $\bdd E'_0$ is a sphere bounding a ball $B$ whose interior is disjoint from $E$.  If no component of $F_i$ lies in $B$, $E_0$ can be isotoped past $E'_0$, reducing $n$ by at least one and maybe more, thereby completing the inductive step.  If some components of $F_i$ lie in $B$, then follow the recipe 
given in Case 1.  For example, if all components of $F_i \cap B$ eventually bound $3$-handles that lie in $B$,  push $M \cap (B \times  [t, 1])$ vertically down to just above level $t$, do the isotopy, then raise $M \cap (B \times  [t, 1])$ back up again.  

\bigskip

{\bf Case 3:}  $k > 0$

\bigskip

Depending on whether $n = 0$ or $n \geq 1$, let $S$ and $B$ be the reducing sphere and $3$-balls described in cases 1) and 2) above.  The inductive hypothesis and a standard innermost disk argument tells us that any $2$-handle attached above level $t_i$ can, starting from highest to lowest, be replaced by a $2$-handle disjoint from the sphere $S$.  Suppose $t_j$ is the level of the first $3$-handle; after the replacement the entire product $S \times [t, t_j)$ is disjoint from $M$.  In fact, following the argument of Case 2, we can isotope the $3$-handles of $M$ (possibly rearranging the ordering of the $3$-handles) so that all of $S \times [t, 1]$ is disjoint from $M$.  Then push $M \cap (B \times [t, 1])$ down to just above level $t$, do the isotopy across $B$ as described in Case 1) or 2), then precisely restore the height of $M \cap (S^3 \times [t, 1])$.
\end{proof}

Note that since all moves are by isotopy, $X$ and $Y$ don't change.  

\begin{lemma}[Torus Unknotting Lemma] \label{lemma:torreimbed}  
Suppose that $Y^*_0$ (or, symmetrically, $X^*_0$) lies in a knotted solid torus $W \subset S^3_0$.  Let $h_0: W \map S^1 \times D^2$ be an orientation-preserving homeomorphism to the unknotted solid torus $S^1 \times D^2 \subset S^3_0$.  Then there is a reimbedding $h: Y \map S^4$ to a rectified critical level embedding so that $h(Y) \cap S^3_0 = h_0(Y^*_0)$ and the number of handles of each index is unchanged.  For $t$ any generic height between the highest $0$-handle and lowest $3$-handle, both of $M^{\pm}_t$ are unchanged.  
\end{lemma}

\begin{proof}   If $\bdd W$ compresses in $X^*_0$, there is nothing to prove: $\bdd$-reduce $W$ to get a $3$-ball, which can be isotoped into $S^1 \times D^2 \subset S^3_0$ and that same isotopy can be applied at every level of $S^3 \times I \subset S^4$.  

So we assume that $\bdd W$ is incompressible in $X^*_0$.  In each successively increasing critical level $t_i  > 0$ ask whether the $2$-handle attached at $t_i$ can be replaced, as in the Prairie-Dog Lemma \ref{lemma:prairie}, by a $2$-handle that lies in $W$.  If it can be done, then do so.  This may alter the critical level embedding of $M$, but only above level $t_{i-1}$.  If success is possible at the critical levels of all $2$-handles, the same can be accomplished for the $3$-handles, as described in Case 2) of the proof of Lemma \ref{lemma:prairie}.  Similarly, at each successively decreasing critical level $t_j < 0$ try to replace cocores of $1$-handles by disks that lie in $W$.  If this can be done for all $1$-handles, then also replace cocores of $0$-handles by $3$-balls in $W$.  If success is possible for all $1$- and $2$- handles, hence at all levels, then $M \cap (\bdd W \times [-1, 1]) = \emptyset$ and so $Y \subset W \times [-1, 1]$.  Then the function $h_0 \times [-1, 1]$ on $ W \times [-1, 1]$, restricted to $Y \subset W \times [-1, 1]$, is the required reimbedding.

We are left with the case where successful replacement of the core of a $2$-handle or cocore of a $1$-handle is not always possible.  Suppose, without loss of generality, that $t_i  > 0$ is the lowest critical level for which the core of the associated $2$-handle cannot be replaced by one that lies in $W \subset S^3_{t_i}$.  Without loss, we assume that the replacements of lower $2$-handles have been done, so $Y \cap (S^3 \times [0, t_i - \eee]) \subset W \times [0, t_i - \eee]$. In particular, the core of the $2$-handle must lie in $X^*_{ t_i - \eee}$ . 

 Choose a disk $D \subset X^*_{t_i - \eee}$ so that its boundary is the same as that of the core of the $2$-handle and, among all such disks, $D$ intersects $\bdd W$ in as few components as possible.  An innermost circle of $D \cap \bdd W$ then cuts off a subdisk of $D$ whose boundary is essential in $\bdd W$ (by choice of $D$) so, since $W$ is knotted, the subdisk must be  a meridian disk for the solid torus $W$.  So at the generic level $t_i - \eee$, $\bdd W$ compresses in  $X^*_{t_i - \eee} \cap W$.   In particular, $Y^*_{t_i - \eee}$ lies in a $3$-ball $B \subset W$.  It is a classical result that simple coning extends the homeomorphism $h_0|B: B \map h_0(B) \subset S^1 \times D^2 \subset S^3$ to an orientation-preserving homeomorphism $H: S^3 \map S^3$.  Define then the embedding of $Y^+_0$ into $S^3 \times [0, 1]$ so that $h|Y \cap (S^3 \times [0, t_i - \eee]) = h_0 \times [0, t_i - \eee]$ and, for $t \geq t_i - \eee$,  $h|Y^*_t = H|Y^*_t$. 

The same argument applies symmetrically to construct $h|Y^-_0$.  Either all co-cores of $1$-handles can be replaced by disks in $W$, in which case we afterwards simply use $h_0$ at every level $t \in [-1, 0]$ or there is a highest critical level $t_i$ for which the cocore of the associated handle on $M$ cannot be replaced by a disk in $W$ and we apply the symmetric version of the construction above. \end{proof}

Here then is an alternative proof of Corollary \ref{cor:genus2}:

\begin{proof} Like any surface in $S^3$, the genus $2$ surface $M^*_0$ compresses in $S^3$, and so it compresses into either $X^*_0$ or $Y^*_0$, say the former.  Maximally compress $M^*_0$ in $X^*_0$.  If $X^*_0$ is a handlebody then Proposition \ref{prop:midHeeg} says $X \cong B^4$, hence its complement $Y \cong B^4$. 

If $X^*_0$ is not a handlebody, then the surface $F$ resulting from maximally compressing the surface $M^*_0$ into $X^*_0$  consists of one or two tori.  Like any surface in $S^3$, $F$ compresses in $S^3$.  The torus component of $F$ that compresses in the complement of $F$ bounds a solid torus $W$ on the side on which the compressing disk lies.  That side can't lie in $X^*_0$, since $F$ is maximally compressed in that direction, so $W$ must contain $M^*_0$ and so indeed all of $Y^*_0$.  The solid torus $W$ is knotted, else $F$ would still compress further into $X^*_0$.  Now apply the Torus Unknotting Lemma \ref{lemma:torreimbed} to reimbed $Y$ in $S^4$ in a level-preserving way so that afterwards $W$ is unknotted; in particular, afterwards $F$ does compress further into the (new) complement of $Y^*_0$.  After perhaps a further iteration of the argument (when $F$ originally consisted of two tori) we have a level-preserving re-embedding of $Y$ in $S^4$ so that afterwards its complement is a handlebody.  It follows from Proposition \ref{prop:midHeeg} then that after such a reimbedding $S^4 - Y \cong B^4$, hence also $Y \cong B^4$.
\end{proof}

\section{Straightening connecting tubes between tori in $S^3$ }

Enlightened by Corollary \ref{cor:HeegU}, observe that there is no generalized Property R obstacle to applying Proposition \ref{prop:midHeeg} to the proof of the genus $3$ Schoenflies Conjecture.  All that is needed is a sufficiently powerful version of the Torus Unknotting Lemma \ref{lemma:torreimbed} that would instruct us how to reimbed some complementary component $Y$ of a genus three $S^3$ in $S^4$ so that its new complement in $S^3_0$ more closely resembles a handlebody (eg it $\bdd$-reduces to a surface of lower genus than $X^*_0$ or $Y^*_0$ did originally).  

Fuelling excitement in this direction is the classic theorem of Fox \cite{Fo} that any compact connected $3$-dimensional submanifold of $S^3$ can be reimbedded as the complement of handlebodies.  What seems difficult to find is a way to extend such a reimbedding of $Y^*_0$ to all of $Y$, as is done in Lemma \ref{lemma:torreimbed}.   It is crucial in the proof of Lemma \ref{lemma:torreimbed} that a solid torus has a unique meridian, whereas of course higher genus handlebodies have infinitely many meridians.  

For genus $3$ embeddings, there is indeed enough information to make such a reimbedding strategy work.  The key is a genus $2$ analogue of the Torus Unknotting Lemma, called the Tube Straightening Lemma.   We precede it with a preparatory lemma in $3$-manifold topology, whose proof is reminiscent of that in \cite{Ga3} or \cite{Th}:

\begin{lemma}  \label{lemma:tubeprep} Suppose $F \subset S^3$ is a genus two surface and $k \subset F$ is a separating curve in $F$.  Denote the complementary components of $F$ by $U$ and $V$ and suppose $k$ bounds a disk $E$ in $V$ so that $U \cup \eta(E)$ is reducible.    Then either 

\begin{itemize}

\item Any simple closed curve in $\bdd N$ that bounds a disk in $S^3 - N$ bounds a disk in $V-N$

\item $N$ can be isotoped in $V$ to be disjoint from $E$ or

\item $k$ bounds a disk in $U$.

\end{itemize}
\end{lemma}

\begin{proof}  Suppose some component of $\bdd N$ is a sphere $S$.  Since $V$ is irreducible, $S$ bounds a ball $B$ in $V$.  Nothing is lost by adding $B$ to $N$, if $N$ is incident to $S$ on the outside of $B$, or removing $B \cap N$ from $N$, if $N$ is incident to $S$ on the inside of $B$.  So we may as well assume that $\bdd N$ has no sphere components.  Essentially the same argument shows that we may take $V - N$ to be irreducible.  For if $S$ is a reducing sphere bounding a ball $B$ in $V$, then any curve in $B \cap \bdd N$ that bounds a disk in $S^3 - N$ bounds a disk in $B - N \subset V - N$ and any curve in $\bdd N - B$ that bounds a disk in $(V - N) \cup B$ also bounds a disk in $V - N$, so without loss, we may delete $B \cap N$ from $N$.   The proof then will be by induction on $-\chi(\bdd N) \geq 0$, assuming now that $V - N$ is irreducible.  

\bigskip

{\bf Case 1:}  $\bdd N$ does not compress in $V - N$

\bigskip

If $\bdd N$ doesn't compress in $S^3 - N$ either, then the first conclusion holds vacuously.  Suppose then that  $S^3 - N$ is $\bdd$ reducible.  A reducing sphere $S$ for $U \cup \eta(E)$ must separate the two tori that are obtained from $F$ by compressing along $E$, since otherwise one ball that $S$ bounds in $S^3$ would lie entirely in $U \cup \eta(E)$.  This implies that $S$ intersects $V$ in an odd number of copies of $E$ or, put another way, it intersects $U$ in a (perhaps disconnected) planar surface with an odd number of boundary components lying on a regular neighborhood of $k$ in $F$.  Let $P$ be a component of the planar surface that has an odd number of boundary components on $F$.  Consider then the result of $0$-framed surgery on $k$ in the manifold $S^3 - N$: $P$ can be capped off to give a sphere which is non-separating in the new manifold, since a meridian of $\eta(k)$ intersects $P$ in an odd number of points.  On the other hand,  $S^3 - N$ itself is $\bdd$-reducible.  No options a)-e) in \cite[Theorem 6.2]{Sch2} are consistent with this outcome (in particular the manifold called $M'$ there having a non-separating sphere)  so we conclude that $S^3 - (N \cup \eta(k))$ is either reducible or $\bdd$-reducible.  In the latter case, consider how a $\bdd$-reducing disk $D$ would intersect the surface $F - \eta(k)$.   We know the framing of $F \cap \bdd \eta(K)$ is a $0$-framing (since each component of $F - \eta(k)$ is a Seifert surface for $k$) so $D \cap F$, if non-empty, consists entirely of simple closed curves.  An innermost disk in $D$ cut off by the intersection (perhaps all of $D$) lies either in $U$ or $V - N$.  But each component of $F - \eta(K)$ is a once-punctured torus, so if it compresses in $U$ or $V-N$ so does its boundary, ie a copy of $k$.  Hence we have that $k$ bounds a disk in either $U$ or $V - N$.  Since $V$ is irreducible,  in the latter case the disk can be isotoped to $E$ in $V$, thereby isotoping $N$ in $V$ off of $E$.  

The same argument applies if $S^3 - (N \cup \eta(k))$ is reducible:  since both $U$ and $V - N$ are irreducible, such a reducing sphere must intersect $F - \eta(k)$ and an innermost disk cut off by the intersection leads to the same contradiction.  

\bigskip

{\bf Case 2:}  $\bdd N$ compresses in $V - N$

\bigskip

The proof is by contradiction:  Choose an essential curve $\ell \subset \bdd N$ and a compressing disk $D$ for $\bdd N$ in $V - N$ so that $\ell$ bounds a disk $D_{\ell}$ in $S^3 - N$ but does not bound a disk  in $V-N$ and, among all such choices of $\ell, D_{\ell}, D$, $|D \cap D_{\ell}|$ is minimal.  If $D$ and $D_{\ell}$ are disjoint, then let $N' = N \cup \eta(D)$.  If there is then a sphere component of $\bdd N'$, the ball it bounds in $V$ can, without loss, be deleted from $N'$.  In any case (perhaps after deleting the ball if a sphere appears in $\bdd N'$), $-\chi(\bdd N') < - \chi(\bdd N)$ and the inductive hypothesis holds for $N'$.  But the conclusion for $N'$ implies the conclusion for $N$, which is contained in $N'$ (eg if $\ell$ bounds a disk in $V - N'$ it bounds a disk in $V - N$), so this is impossible.

If $D$ and $D_{\ell}$ are not disjoint, note that all curves of intersection must be arcs, else an innermost disk cut off in $D$ could be used to surger $D_{\ell}$ and would lower $|D \cap D_{\ell}|$.  Similarly, let $D'$ denote the disk cut off from $D$ by an outermost arc of $D \cap D_{\ell}$ in $D$ and  let $D_{\ell'}, D_{\ell''} \subset V - N$ denote the two disks obtained by the $\bdd$-compression of $D_{\ell}$ to $\bdd N$ along $D'$.  Both of these disks intersect $D$ in fewer components than $D_{\ell}$ did, so by choice of $\ell$ and $D_{\ell}$ the boundary of each new disk must bound a disk in $V-N$ (of course, if either curve is inessential in $\bdd N$ then it automatically bounds a disk in $\bdd N \subset V$).   A standard innermost disk argument shows that $D_{\ell'}$ and $D_{\ell''}$  can be taken to be disjoint; band them together via the band in $\bdd N$ which undoes the $\bdd$-compression by $D'$.  To be more explicit, note that the arc $\bdd D' \cap \bdd N$ defined a band move on $\bdd D_{\ell}$ that split $\bdd D_{\ell}$ into $\bdd D_{\ell'}$ and $\bdd D_{\ell''}$.  Undo that band move to recover the curve $\bdd D_{\ell}$, now bounding a disk (namely the band sum of $D_{\ell'}$ and $D_{\ell''}$) that lies in $V - N$. This contradicts our original choice of $\ell$.
\end{proof}

\begin{lemma}[Tube Straightening Lemma] \label{lemma:tubereimbed}  
Suppose that $Y \cap S^3_0$ (or, symmetrically, $X \cap S^3_0$) lies in $V \subset S^3_0$, with closed complement $U$, and $\bdd U = \bdd V$ is of genus $2$.   Suppose $V$ contains a separating compressing disk $E$ so that the manifold $U_+$ obtained from $U$ by attaching a $2$-handle along $E$ is reducible.  Then there is an embedding $h_0: V  \map S^3$ so that $h_0(\bdd E)$ bounds a disk in $U'$ the complement of $h_0(V)$.  There is also a reimbedding $h: Y \map S^4$ to a rectified critical level embedding so that 
\begin{itemize}
\item $h(Y) \cap S^3_0 = h_0(Y^*_0)$
\item the number of handles of each index is unchanged and
\item for $t$ any generic height between the highest $0$-handle and lowest $3$-handle, both of $M^{\pm}_t$ are unchanged.    
\end{itemize}
\end{lemma}

\begin{proof}   Let $h_0: V \map S^3$ be the reimbedding (unique up to isotopy) that replaces the $1$-handle in $V$ that is dual to $E$ with a handle intersecting the reducing sphere for $U_+$ in a single point. Then after the reimbedding $\bdd E$ bounds a disk in the complement of $h_0(V)$,  namely the complement of $E$ in the reducing sphere.  

 In each successively increasing critical level $t_i  > 0$ ask whether the $2$-handle attached at $t_i$ can be replaced, as in the Prairie-Dog Lemma \ref{lemma:prairie}, by a $2$-handle that lies in $V$.  If it can, then do so.  This may alter the critical level embedding of $M$, but only above level $t_{i-1}$.  If success is possible at the critical levels of all $2$-handles, the same can be accomplished for the $3$-handles, as described in Case 2) of the proof of Lemma \ref{lemma:prairie}.  Similarly, at each successively decreasing critical level $t_j < 0$ try to replace cocores of $1$-handles by disks that lie in $V$.  If this can be done for all $1$-handles, then also replace cocores of $0$-handles by $3$-balls in $V$.  If success is possible for all $1$- and $2$- handles, hence at all levels, then $M \cap (\bdd V \times [-1, 1]) = \emptyset$ and so $Y \subset V \times [-1, 1]$.  Then the function $h_0 \times [-1, 1]$ on $ V \times [-1, 1]$, restricted to $Y \subset V \times [-1, 1]$, is the required reimbedding.

We are left with the case where successful replacement of the core of a $2$-handle or cocore of a $1$-handle is not always possible.  Suppose, without loss of generality, that $t_i  > 0$ is the lowest critical level for which the core of the associated $2$-handle cannot be replaced by one that lies in $V \subset S^3_{t_i}$.  Without loss, we assume that the replacements of lower $2$-handles have been done, so $Y \cap (S^3 \times [0, t_i - \eee]) \subset (V \times [0, t_i - \eee])$. In particular, the core of the $2$-handle must lie in $X^*_{ t_i - \eee}$ . 

Now apply Lemma \ref{lemma:tubeprep} using $Y^*_{t_i - \eee}$ for $N$.  By assumption, the boundary of the core of the $2$-handle bounds no disk in $V \cap X^*_{ t_i - \eee}$ so the first possible conclusion of Lemma \ref{lemma:tubeprep} cannot hold.  If the last holds, there was nothing to prove to begin with.  (Take $h = identity$.) Hence we conclude that the second conclusion holds: $Y^*_{t_i - \eee}$ can be isotoped to be disjoint from $E$.  But once this is true, the reimbedding $h_0$ has no effect on $Y^*_{t_i - \eee}$; that is, $Y^*_{t_i - \eee}$ is isotopic to $h_0(Y^*_{t_i - \eee})$.  Hence we can define $h|(Y \cap (S^3 \times [0, 1]))$ to be $h_0 \times [0, {t_i - \eee}]$ on $Y \cap (S^3 \times [0, {t_i - \eee}])$, followed by a quick isotopy of $h_0(Y^*_{t_i - \eee})$ to $Y^*_{t_i - \eee/2}$ followed by the unaltered embedding above level $t_i - \eee/2$.  Note that this unaltered embedding is not necessarily the original embedding, because of changes made while ensuring that earlier $2$-handles lie in $V$.  

Finally, the argument can be applied symmetrically on $Y \cap (S^3 \times [-1, 0]).$   \end{proof}

\section{Weak Fox reimbedding via Unknotting and Straightening}

In this section we show that, for a genus $3$ surface in $S^3$, the operations of Torus Unknotting and Tube Straightening described above suffice to give a weak version of Fox reimbedding.  That is, for a genus $3$ surface $F \subset S^3$ there is a sequence of such reimbeddings, not necessarily all operating on the same complementary component of $F$, so that eventually a complementary component is a handlebody.  Although the context of this section appears to be $3$-manifold theory, the notation is meant to be suggestive of the eventual application to the genus $3$ Schoenflies Conjecture.  In particular, the terms {\em torus unknotting} and {\em tube straightening} as used in this section refer to the $3$-dimensional reimbedding $h_0$ given in, respectively, Lemma \ref{lemma:torreimbed} and Lemma \ref{lemma:tubereimbed}.

Suppose $F \subset S^3$ is a surface dividing $S^3$ into two components denoted $X$ and $Y$. 
Suppose $D_1$ is a compressing disk for $F$ in $X$ giving rise to a new surface $F_1 \subset S^3$ with complementary components $X_1 = X - \eta(D_1)$ and $Y_1 = X \cup \eta(D_1)$.  Suppose $D_2$ is a compressing disk for $F_1$ in $X_1$ or $Y_1$, giving rise to a new surface  $F_2 \subset S^3$ with complementary components $X_2$ and $Y_2$. Continue to make such a series of compressions via compressing disks $D_i, i = 1, ... , n$ so that each $D_i$ lies either in $X_{i-1}$ or $Y_{i-1}$ (the complementary components of $F_{i-1} \subset S^3$) until all components of $F_n$ are spheres.  

\begin{defin}  $F$ {\em can be straightened} if there is a sequence of torus unknottings and tube straightenings of $F$, as described in Lemmas \ref{lemma:torreimbed} and Lemma \ref{lemma:tubereimbed}, so that afterwards either 
\begin{itemize}
\item  $X$ or $Y$ is a handlebody or
\item the order of the compressions $D_1,..., D_n$ can be rearranged so that afterwards some $\bdd D_i$ is inessential in $F_{i-1}$ and so can be eliminated from the sequence.
\end{itemize}
In the spirit of Lemma \ref{lemma:prairie}, any disk $D_i$ can be replaced by a disk $D'_i$ in the same complementary component, so long as $\bdd D_i$ and $ \bdd D'_i$ are isotopic  in $F_{i-1}$.  
\end{defin}

The following series of lemmas assumes we are given such a sequence for $F \subset S^3$, with the first compressing dis $D_1$ in $X$.  There are, of course, symmetric statements if $D_1 \subset Y$.

\begin{lemma}  \label{lemma:sameside}   If $D_2 \subset X_1$ then  $F$ can be straightened.  More specifically, there is a sequence of torus unknottings which convert $X$ into a handlebody \end{lemma}

\begin{proof}  Since $X_1 = X - \eta(D_1)$, the disks $D_1$ and $D_2$ can be thought of as disjoint disks in $X$ and the compressions given by $D_1, D_2$ can be performed simultaneously.  $F_2$ consists of one, two, or three tori, depending on how many of $D_1, D_2$ are non-separating. The proof is by induction on the number of $D_1, D_2$ that are separating disks.  

If there are none, so both $D_1, D_2$ are non-separating, then $F_2$ is a single torus; if the solid torus it bounds lies in $X_2$ then the original $X$ was a handlebody and we are done.  If the solid torus that $F_2$ bounds lies in $Y_2$ then all of $F$ also lies in that solid torus.  After a torus unknotting, $F_2$ also bounds a solid torus in $X_2$ and $X$ is a handlebody, as required.  

If $D_1$ (or $D_2$) is separating then one of the components $T$ of $F_1$ is a torus.  If $T$ bounds a solid torus in $X_1$ then we may as well have used the meridian of that torus for $D_1$ and invoked the inductive hypothesis.  If, on the other hand, $T$ bounds a solid torus in $Y_1$ then all of $F$ also lies in that solid torus.  After perhaps a torus unknotting, $T$ bounds a solid torus in $X_1$ as well, and again we could replace $D_1$ by a meridian of that torus and invoke the inductive hypothesis.  \end{proof}

\begin{lemma}  \label{lemma:oppside}   If $D_1 \subset X$ and $D_2 \subset Y_1$ are both separating, and $D_2 \subset Y \subsetneq Y_1$ (ie $D_2$ does not pass through the $1$-handle dual to $D_1$) then  $F$ can be straightened.  \end{lemma}

\begin{proof}  Since the interiors of both $D_1, D_2$ are disjoint from $F$, their order can be rearranged, so there is symmetry between the two.  Since $D_1$ is separating there is a torus component $T_1 \subset F_1$ and $T_1$ bounds a component $W_1$ of $X_1$ whose interior is disjoint from $F_2$.  After perhaps a torus unknotting of its complement we may as well assume that $W_1$ is a solid torus.  Eventually some compression $D_i$ will compress $T_1$ to a sphere; if $D_i \subset W_1$ then we could have done $D_i$ before $D_1$, making $D_1$ (coplanar to $D_i$) redundant, and thereby reduced the number of compressions. Thus we may as well assume $D_i \subset Y_{i-1}$, so $W_1$ is an unknotted solid torus.  Similarly, the torus component $T_2 \subset F_2$ not incident to $D_1$ bounds a solid unknotted torus $W_2 \subset Y_2$.  See Figure \ref{fig:oppside}.

 \begin{figure}[tbh]
  \centering
  \includegraphics[width=.7\textwidth]{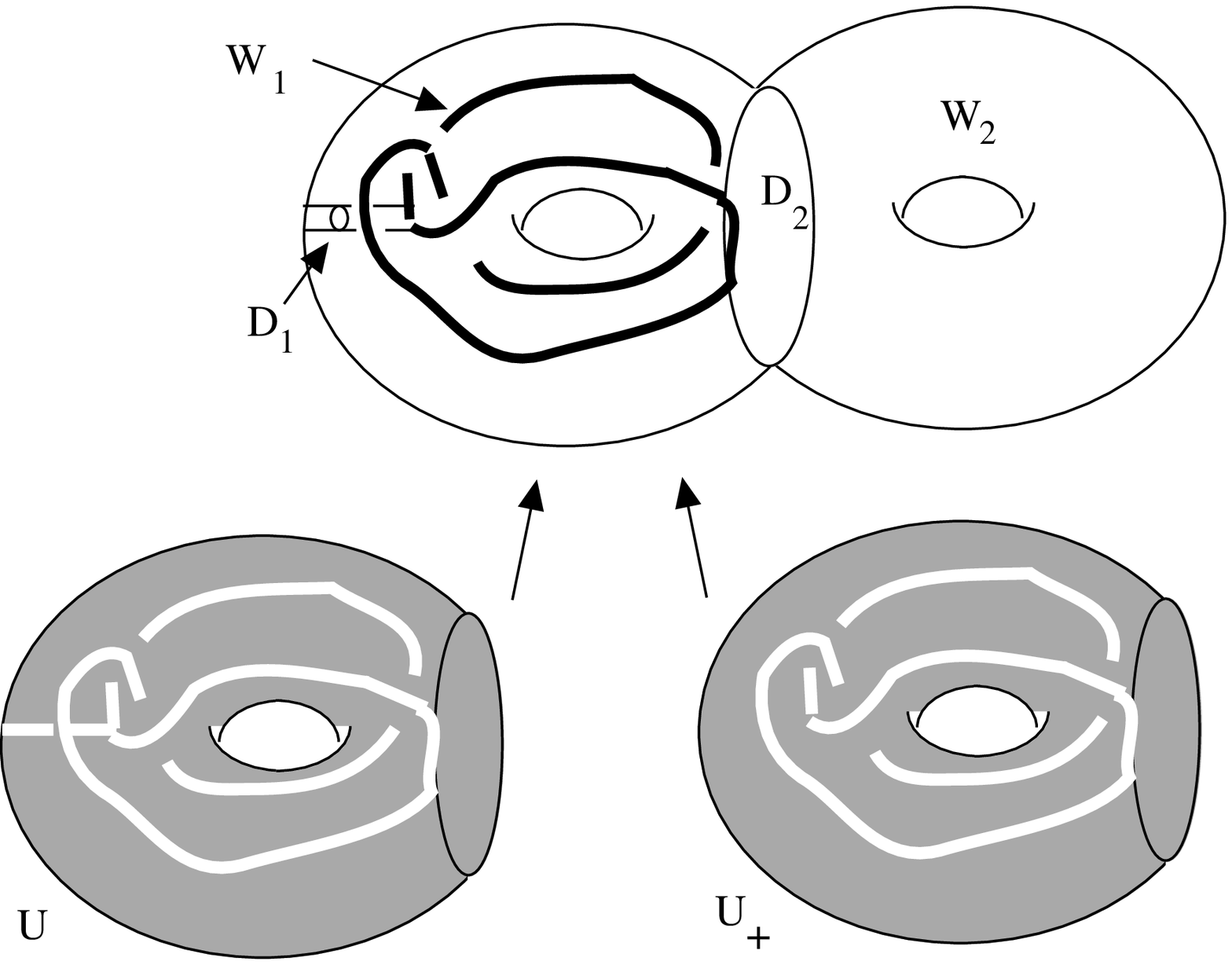}
  \caption{} \label{fig:oppside}
  \end{figure}  

Assume, with no loss of generality, that $D_3 \subset Y_2$.  We have already seen that $D_3$ can't compress $T_2$, else the number of compressions could be reduced.  It follows that $D_3$ either compresses $T_1$ or it compresses the third torus $T_3$ created from $F$ by the compressions $D_1, D_2$.  Either case implies that the component $U_+ = Y_2 - W_2$ of $Y_2$ is reducible.  Apply tube straightening to $U = U_+ \cap Y$ so that the tube dual to $D_1$ passes through the reducing sphere of $U_+$ once.  After the straightening, the disk $D_i$ that eventually compresses $T_1$ into $Y_i$ lies in $Y$, so it can be the first compression, making $D_1$ redundant.
 \end{proof}
 
  \begin{lemma}  \label{lemma:oppside2}   If $D_1 \subset X$ and $D_2 \subset Y_1$ are both separating and if $D_3$ or both $D_4$ and $D_5$ are towards the $X$ side, then  $F$ can be straightened.
\end{lemma}

\begin{proof}  In this case there is not necessarily symmetry between $D_1$ and $D_2$, but the argument of Lemma \ref{lemma:oppside} still applies if the compression disk for $T_3$ lies on the $X$- side rather than the $Y$-side or if $T_2$ compresses before $T_3$.  So that is what we now verify:  If $D_3$ lies on the $X$-side it  compresses $T_2$ or $T_3$ into $X_2$ and we are done.  On the other hand, if $D_4$ and $D_5$ both lie on the $X$-side then since $T_1$ compresses on the $Y$-side, one of $D_4$ or $D_5$ is the compression disk for $T_3$. \end{proof}

 \begin{lemma}  \label{lemma:oppside3}   If $D_1 \subset X$ is separating, and $Y$ is $\bdd$-reducible, then  either $D_2 \subset Y \subsetneq Y_1$ or $F$ can be straightened.
\end{lemma}

\begin{proof} Consider the torus component $T_1$ of $F_1$ and the component $W_1 \subset X_1$ it bounds.  As noted above, the interior of $W_1$ is disjoint from $F_1$ and so, perhaps after a torus unknotting, we can assume that $W_1$ is a solid torus.  Moreover, the disk $D_i$ that eventually compresses $T_1$ lies on the $Y$-side and not in $W_1$.  Hence $W_1$ is an unknotted solid torus.

Following Lemma \ref{lemma:sameside} we may as well assume that $D_2 \subset Y_1$.  By assumption $Y$ is $\bdd$-reducible and, after attaching a $2$-handle (a neighborhood of $D_1$) to get $Y_1$ the resulting manifold still is $\bdd$-reducible, via $D_2$.  It follows from the Jaco handle addition lemma \cite{Ja} that there is a $\bdd$-reducing disk $J \subset Y$ for $Y$ whose boundary is disjoint from $\bdd D_1$.  Take $J$ to be non-separating, if this is possible.  

If $\bdd J$ lies on $T_1$ then it is parallel to the disk $D_i$ that eventually compresses $T_1$ so we may as well do that compression before $D_1$, making $D_1$ redundant and so reducing the number of compressions.  So we henceforth suppose $\bdd J$ lies on the other, genus $2$ component of $F - \bdd D_1$, and so lies on $\bdd F_1 - T_1$.  

If $\bdd J$ is inessential in $F_1$ then it is parallel to $\bdd D_1$ in $F$.  Thus $\bdd D_1$ also bounds the disk $J \subset Y$ and the union of the two gives a reducing sphere for $Y_1$ that intersects the $1$-handle dual to $D_1$ in a single point.  It follows that the longitude $\bdd D_i$ of $T_1$ bounds a disk in $Y$ and we are done as before.

If $\bdd J$ is essential in $F_1$ compress $Y$ along $J$ to get $Y_J$ and consider the component $W_J \subset Y_J$ (in fact all of $Y_J$ if $J$ is non-separating) such that $\bdd D_1 \subset \bdd W_J$ and $\bdd W_J$ has genus $2$.   The manifold $W^+_J$ obtained by attaching $\eta(D_1)$ to $W_J$ has boundary the union of two tori.  $W^+_J$ can also be viewed as a component of $Y_1 - \eta(J)$.  

We claim that either $D_2 \subset Y$ or $W^+_J$ is reducible.  This is obvious if there is a disk component of $D_2 - J$ that can't be removed by an isotopy or if $\bdd D_2 \cap J = \emptyset$.  The alternative then is that $D_2 \cap J$ is a non-empty collection of arcs.  Consider an outermost disk $D'$ cut off from $D_2$ by $J$.  We may as well assume $D' \subset W^+_J$.  For if it's not then $J$ is separating, $D'$ is a meridian of the other complementary component (a solid torus) and we may as well have used $D'$ for $J$, thereby eliminating the case that $J$ is separating. With then $D' \subset W^+_J$, compress the torus boundary component of $\bdd W^+_J$ along $D'$ to get a reducing sphere for $W^+_J$.  

Finally, if $W^+_J$ is reducible, apply tube straightening to the surface $\bdd W_J$, replacing the handle dual to $D_1$ by a handle intersecting the reducing sphere once, allowing the same conclusion as before.   \end{proof}

 \begin{cor}  \label{cor:oppside1}   If $D_1 \subset X$ and $D_2 \subset Y_1$ are both separating, and $Y$ is $\bdd$-reducible, then  $F$ can be straightened.
\end{cor}

\begin{proof} Combine Lemmas \ref{lemma:oppside3} and \ref{lemma:oppside}. \end{proof}



\begin{cor} \label{cor:separate}   If $D_1 \subset X$ is separating then either:

\begin{itemize}

\item at least two of the three non-separating compressing disks are on the $Y$-side and, if $D_2$ is separating, $Y$ is $\bdd$-irreducible or

\item  $F$ can be straightened.

\end{itemize}
\end{cor}

\begin{proof}  Following Lemmas \ref{lemma:sameside}, \ref{lemma:oppside2}, and \ref{lemma:oppside3} we may as well assume that $D_2 \subset Y \subsetneq Y_1$ and $D_2$ is non-separating.  As before, let $T_1$ be the torus component of $F_1$ bounding a component $W_1$ of $X_1$ whose interior is disjoint from $F_1$.  

If the compressing disk $D_i$ that eventually compresses $T_1$ lies on the $X$-side, then $D_i \subset W_1$ (and $W_1$ is a solid torus).  No earlier compressing disk can be incident to $T_1$ so in fact $D_i$ could have been done before $D_1$, making the latter redundant.  This reduces the number of compressions.  

If, on the other hand, $D_i$ lies on the $Y$-side then both $D_i$ and $D_2$ are non-separating disks lying on the $Y$-side, as required.
\end{proof}

\begin{lemma}  \label{lemma:nonsep}   If $D_1 \subset X$ is non-separating, then  either all succeeding disks $D_2, D_3$ are non-separating or $F$ can be straightened.
\end{lemma}

\begin{proof}  Since $D_1$ is non-separating, $F_1$ is a genus $2$ surface.  Its complementary component $Y_1$ contains all of $Y$.  If $D_2$ is also non-separating then $F_3$ is a torus, for which any compressing disk is non-separating, giving the result.

So suppose $D_2$ is separating.  Following Lemma \ref{lemma:sameside} we may as well assume $D_2 \subset Y_1$, so $\bdd Y_2$ consists of two tori.  Hence $D_3$ is non-separating and compresses one of the tori.  If $D_3 \subset Y_2$ then $D_3$ could have been done before $D_2$, making $D_2$ redundant.  If $D_3 \subset X_2$ then the result of the compression is a sphere which could have been viewed as a reducing sphere for $Y_2$.  Apply tube straightening, using $Y_1$ for $N$.  After that reimbedding, $D_3 \subset X_1 \subsetneq X_2$ does not pass through $D_2$ so the compression along $D_3$ could be done before the compression along $D_2$, making $D_2$ redundant. \end{proof}

We note in passing, though the fact will not be used, that if there are non-separating compressing disks on both sides, they may be taken to be disjoint:

 \begin{prop}  \label{prop:oppside4}   If $D_1 \subset X$ is non-separating, and $Y$ has a non-separating $\bdd$-reducing disk $E$, then  either $\bdd D_1 \cap \bdd E = \emptyset$ or $F$ can be straightened.
\end{prop}

\begin{proof} Following Lemma \ref{lemma:sameside} we may as well assume that $D_2 \subset Y_1$.  By assumption $Y$ is $\bdd$-reducible and, after attaching a $2$-handle (a neighborhood of $D_1$) to get $Y_1$ the resulting manifold still $\bdd$-reducible, via $D_2$.  It follows from the Jaco handle addition lemma  \cite{Ja} that there is a $\bdd$-reducing disk $J \subset Y$ for $Y$ whose boundary is disjoint from $\bdd D_1$.  Take $J$ to be non-separating, if this is possible.  

Suppose first that $\bdd J$ is inessential in $F_1$.  Then the disk it bounds in $F_1$ contains both copies of $D_1$ resulting from the compression of $F$ along $D_1$.  Put another way, $J$ cuts off a component $W_J$ from $Y$ that has torus boundary and whose interior is disjoint from $F$.  Following, perhaps, a torus unknotting, we may assume that $W_J$ is a solid torus.  A standard innermost disk, outermost arc argument between $J$ and $D_2$ ensures that they can be taken to be disjoint, so $D_2$ compresses the other, genus $2$ component of $Y - \eta(J)$.  This compression, together with the compression via the meridian of $W_J$, $\bdd$-reduces $Y$ to one or two components, each with a torus boundary.  After perhaps some further torus unknottings, $Y$ becomes then a handlebody, as required. 

So suppose henceforth that $\bdd J$ is essential in $F_1$.   Compress $Y$ along $J$ to get $Y_J$ and consider the component $W_J \subset Y_J$ (in fact all of $Y_J$ if $J$ is non-separating) such that $\bdd D_1 \subset \bdd W_J$ and $\bdd W_J$ has genus $2$.   The manifold $W^+_J$ obtained by attaching $\eta(D_1)$ to $W_J$ has boundary a torus.  $W^+_J$ can also be viewed as a component of $Y_1 - \eta(J)$.  

Consider an outermost disk $E'$ of $E$ cut off by $J$, or let $E' = E$ if $E$ is disjoint from $J$.  We may as well assume that $E'$ lies in $W_J$, since if $J$ is separating and $E'$ lies in the other component, we should have taken $E'$ for $J$.  If $E'$ is inessential in $W_J$ then $E' = E$ is parallel to $J$ (since $E$ is non-separating), so $\bdd D_1 \cap \bdd E = \emptyset$ and we are done.  If $E'$ is essential in $W_J$ then each component of $Y - (\eta(J) \cup \eta(E'))$ has interior disjoint from $F$ and is bounded by a torus.  Following some torus unknottings we can take them to be solid tori.  In that case $Y$ is a handlebody, as required.  \end{proof} 

\section{The genus $3$ Schoenflies Conjecture}

We now apply the results of the previous sections to complete the proof of the genus $3$ Schoenflies Conjecture.   

\begin{thm} \label{thm:main}  Suppose $M$ is a genus $3$ rectified critical level embedding of $S^3$ in $S^4$.  Then after a series of reimbeddings of $M$ via other genus $3$ rectified critical level embeddings, each of which changes at most one of the complementary components of $M$, one of those complementary components is $B^4$
%




\end{thm}

\begin{proof} We assume that any possible genus $3$ rectified critical level reimbedding of $M$ that preserves at least one complementary component and simultaneously decreases the number of handles has been done.  If any further such a sequence of reimbeddings via Lemmas \ref{lemma:tubereimbed} or \ref{lemma:torreimbed} (reimbeddings that don't raise the number of handles) results in $X^*_0$ (resp  $Y^*_0$) becoming a handlebody, then $X$ (resp $Y$) is $B^4$, via Proposition \ref{prop:midHeeg}.  So we assume no such further sequence exists and use the results of the previous section to see if there are other options.  With no loss of generality, assume the cocore $D_1$ of the highest $1$-handle lies in $X^*_0$ and let $E_1$ denote the core of the lowest $2$-handle.

Following Lemma \ref{lemma:sameside} we can assume that the cocore $D_2$ of the previous $1$-handle lies on the $Y$-side and the cocore $E_2$ of the next $1$-handle lies on the side opposite the one on which $E_1$ lies.  

{\bf Claim:} Perhaps after rearranging the ordering of the handles, at least one of $D_1$ and $E_1$ is non-separating.

Suppose $D_1$ and $E_1$ are both separating.  If both $D_1$ and $E_1$ lie in $X^*_0$, then it follows from Corollary \ref{cor:separate} that at least two of the non-separating cores of the $2$-handles and at least two of the non-separating cocores of the $1$-handles all lie on the $Y$-side of the surfaces to which they are attached.  So both $Y^+_0$ and $Y^-_0$ are $4$-dimensional handlebodies of genus at least $2$.  But the Mayer-Vietoris sequence for $Y^+_0, Y^-_0$ glued along $Y^*_0$ then contradicts $H_*(Y) \cong H_*(B^4)$, cf the proof of Lemma \ref{lemma:rhosum}.

On the other hand, if $D_1 \subset X^*_0$ and $E_1 \subset Y^*_0$ then, following Corollary \ref{cor:oppside1} and Lemma \ref{lemma:oppside3}, $D_2$ is non-separating and lies in $Y^*_0$.  Then interchange $D_1$ and $D_2$, using the $1$-handle dual to $D_2$ as the highest $1$-handle.  The new arrangement establishes the claim.

Following the claim, we can, with no loss, assume that $E_1$ is non-separating. 
Then according to Lemma \ref{lemma:nonsep}, all of the cores of $2$-handles are non-separating, so each surface $F_i$ at or above height $t = 0$ are connected.  Hence there is at most one $3$-handle in $M$ and, passing this $3$-handle over the north pole if necessary, this guarantees that each of $X$ and $Y$ have induced handle structures without $3$-handles.   Following the comments preceding Lemma \ref{lemma:dualrule}, the sum of the genera of the $4$-dimensional handlebodies  $X^-_0$ and $Y^-_0$ is $3$, so one of them, say $X^-_0$, has genus $\leq 1$.  $X$ is then obtained from the genus $0$ or $1$ handlebody $X^-_0$ by attaching some $2$-handles but no $3$-handles.  Moreover, $\bdd X$ is a sphere.  If $genus(X^-_0) = 0$ (resp $genus(X^-_0) =1$) then no (resp $1$) $2$- handle must be attached, to ensure that $H_1(X) = H_2(X) = 0$.  In the first case, since no $2$-handles (hence no handles at all) are attached, $X \cong X^-_0 \cong B^4$.  In the second case, one $2$-handle is attached to $X^-_0 \cong S^1 \times D^3$ and so $X \cong B^4$, by Corollary \ref{cor:PropR}.
\end{proof}

\begin{cor}  \label{cor:final} Each complementary component of a genus three rectified critical level embedding of $S^3$ in $S^4$ is a $4$-ball. \end{cor}

\begin{proof}  Let $M$ be a genus three rectified critical level embedding of $S^3$.  Following Theorem \ref{thm:main} there is a sequence of such embeddings $$M = M_0, M_1, M_2, ... , M_n$$ so that one of the complementary components of $M_n$ is $B^4$ and for each $i = 0,...,n-1$, one of the complementary components of $M_{i}$ is homeomorphic to a complementary component of $M_{i+1}$.  The argument is by induction on $n$, exploiting the fact that the complement of $B^4$ in $S^4$ is $B^4$.  

Since one complementary component of $M_n$ is $B^4$, both complementary components are. If $n = 0$ we are done.  For $n \geq 1$, note that since one complementary component of $M_{n-1}$ is homeomorphic to a complementary component of $M_n$, that complementary component is $B^4$.   This completes the inductive step.
\end{proof}

\end{document}